\documentclass[11pt]{amsart}
\usepackage[arrow, matrix]{xy}
\usepackage{amsmath}
\usepackage{amssymb}
\usepackage{amscd}
\newtheorem{dfn}{Definition}[section]
\newtheorem{rem}[dfn]{Remark}

\newtheorem{thm}[dfn]{Theorem}
\newtheorem{lem}[dfn]{Lemma}
\newtheorem{prop}[dfn]{Proposition}
\newtheorem{cor}[dfn]{Corollary}

\newtheorem{defn}[dfn]{Definition}

\def\no{\noindent}

\def\0{\emptyset}
\def\half{\frac{1}{2}}
\def\proof{\par\medskip\noindent{\it Proof: }}

\def\>{\rangle}
\def\<{\langle}

\def\C{\mathbb C}
\def\t{\widetilde}
\def\R{\mathbb R}

\def\E{{\mathbb E}}

\def\h{\mathbb H}

\def\O{{\mathcal O}}

\def\M{{\mathcal M}}
\def\eps{\epsilon}

\def\ka{\kappa}

\def\B{\mathcal{B}}
\def\mh{\mathcal{H}}
\def\le{\lim_{\eps\to0}}
\def\ha{\widehat}

\def\th{\theta}
\def\al{\alpha}
\def\be{\beta}
\def\ga{\gamma}

\def\del{\delta}

\def\si{\sigma}

\def\r{{\mathcal R}}
\def\la{\lambda}
\def\La{\Lambda}
\def\8{\infty}

\def\om{\omega}
\def\D{\partial}

\def\k{{\mathfrak k}}
\def\b{{\mathfrak b}}
\def\g{{\mathfrak g}}
\def\e{{\mathcal E}}
\def\U{{\mathcal U}}

\def\1{\hbox{\bf 1}}

\def\nab{\nabla}

\begin{document}
\title{The Symplectic
       Geometry of Polygons in Hyperbolic 3-space}
\author[Michael Kapovich]{Michael Kapovich*}
\address{Department of Mathematics\\
     University of Utah\\
     Salt Lake City, UT 84112, USA}
\email{kapovich@math.utah.edu}
\thanks{*\ Research partially supported by NSF grant DMS-96-26633.}
\author[John J. Millson]{John J. Millson$^\dag$}
\address{Department of Mathematics\\
     University of Maryland\\
     College Park, MD 20742, USA}
\email{jjm@math.umd.edu}
\thanks{$^\dag$\ Research partially supported by NSF grant DMS-98-03520.}
\author[Thomas Treloar]{Thomas Treloar}
\address{Department of Mathematics\\
     University of Maryland\\
     College Park, MD 20742, USA}
\email{txt@math.umd.edu}
\date{\today}

\begin{abstract}
We study the symplectic geometry of the moduli spaces $M_r=M_r(\h^3)$ of closed
n-gons
with fixed side-lengths in hyperbolic three-space.  We prove that these moduli
spaces have almost canonical symplectic structures.  They are the symplectic
quotients of $B^n$ by the dressing action of $SU(2)$ (here $B$ is the
subgroup of the Borel subgroup of $SL_2(\C)$ defined
below).  We show that the hyperbolic
Gauss map sets up a real analytic isomorphism between the spaces $M_r$ and the
weighted quotients of $(S^2)^n$ by $PSL_2(\C)$ studied by Deligne and Mostow.
We construct an integrable Hamiltonian system on $M_r$ by bending polygons
along nonintersecting diagonals.  We describe angle variables and the momentum
polyhedron for this system.  The results of this paper are the analogues for
hyperbolic space of the results of \cite{KM2} for $M_r(\E^3)$, the space of
n-gons with fixed side-lengths in $\E^3$.  We prove $M_r(\h^3)$ and
$M_r(\E^3)$ are symplectomorphic.
\end{abstract}

\maketitle

\section{Introduction}
An (open) n-gon $P$ in hyperbolic space $\h^3$ is an ordered
(n+1)-tuple $(x_1,...,x_{n+1})$ of points in $\h^3$ called the
vertices.  We join the vertex $x_i$ to the vertex $x_{i+1}$ by the
unique geodesic segment $e_i$, called the $i$-th edge.  We let
$Pol_n$ denote the space of n-gons in $\h^3$.  An n-gon is said to
be closed if $x_{n+1} = x_1$.  We let $CPol_n$ denote the space of
closed n-gons.  Two n-gons $P=[x_1,...,x_{n+1}]$ and
$P'=[x'_1,...,x'_{n+1}]$ are said to be equivalent if there exists
$g \in PSL_2(\C)$ such that $gx_i = x'_i$, for all $1\leq i \leq
n+1$.  We will either represent an n-gon $P$ by its vertices or
its edges, $P=[x_1,...,x_{n+1}]=(e_1,...e_n).$

Let $r=(r_1,...,r_n)$ be an n-tuple of positive numbers.  This paper is
concerned with the symplectic geometry of the space of closed n-gons in $\h^3$
such that the $i$-th edge $e_i$ has side-length $r_i$, $1 \leq i \leq n$,
modulo $PSL_2(\C)$.  We will assume in this paper (with the exception of
\S\ref{S3}) that $r$ is not on a wall of $D_n$ (see \S\ref{S2}), hence $M_r$
is a real-analytic manifold.

The starting point of this paper is (see \S\ref{S4})

\begin{thm}
The moduli spaces $M_r$ are the symplectic quotients obtained from the dressing
action of $SU(2)$ on $B^n$.
\end{thm}

Here $B=AN$ is the subgroup of the Borel subgroup of $SL_2(\C)$,
$B =\{\left(
\begin{smallmatrix}
\la & z \\
0 & \; \la^{-1}
\end{smallmatrix}
\right) : \la \in \R_+, \; z \in \C\}$.
$B$ is given the Poisson Lie group structure corresponding to the Manin triple
$(\mathfrak{sl}(2,\C),\mathfrak{su}(2),\b)$ with $\< , \>$ on
$\mathfrak{sl}(2,\C)$ given by the imaginary part of the Killing form.

\begin{rem}
As a consequence of Theorem 1.1, the spaces $M_r$ have an almost canonical
symplectic structure (the symplectic structure depends on a choice of Iwasawa
decomposition of $SL_2(\C)$ or a ray in $\h^3$, but given two such choices,
there exist (infinitely many) $g\in SL_2(\C)$ inducing an isomorphism of the two
Poisson structures).
\end{rem}

Our next theorem relates the moduli spaces $M_r$ to the weighted quotients
$Q_{sst}=Q_{sst}(r)$ of $(S^2)^n$ constructed by Deligne and Mostow in
\cite{DM}.  By extending the sides of the n-gon in the positive direction until
they meet $S^2 = \partial_{\8} \h^3$, we obtain a map, the hyperbolic Gauss map
$\gamma : Pol_n \to (S^2)^n$. We then have (here we assume $M_r$ is smooth)

\begin{thm}
The hyperbolic Gauss map induces a real analytic diffeomorphism $\gamma:M_r \to
Q_{sst}(r)$.
\end{thm}

\begin{rem}
In \cite{KM2}, the first two authors constructed an analogous analytic
isomorphism $\ga : M_r(\E^3) \to Q_{sst}(r)$ where $M_r(\E^3)$ is the moduli
space of n-gons with the side-lengths $r=(r_1,...,r_n)$ in Euclidean space
$\E^3$.  Although they gave a direct proof, this latter result was a
consequence of the Kirwan-Kempf-Ness theorem, \cite{Ki},\cite{KN}, relating
Mumford quotients to symplectic quotients.  Our new result (Theorem 1.3 above)
relates a Mumford quotient to a quotient of a symplectic manifold by a
\textbf{Poisson} action.
\end{rem}

The key step (surjectivity) in the proof of Theorem 1.3 is of independent
interest.  We could try to invert $\ga : M_r \to Q_{sst}$ as follows.  Suppose
we are given $\xi=(\xi_1,...,\xi_n) \in Q_{sst}$.  We wish to construct $P \in
M_r$ with $\ga(P)=\xi$.  Choose $x \in \h^3$.  Put the first vertex $x_1=x$.
Let $\si_1$ be the geodesic ray from $x_1$ to $\xi_1$.  Let $x_2$ be the point
on $\si_1$ with $d(x_1,x_2)=r_1$.  Let $\si_2$ be the ray from $x_2$ to $\xi_2$.
Cut off $\si_2$ at $x_3$ so that $d(x_2,x_3)=r_2$.  We continue in this way
until we get
$P=[x_1,...,x_{n+1}]$.  However it may not be the case that $P$ closes up
(i.e. $x_{n+1}=x_1$).

\begin{thm}
Suppose $\xi$ is a stable configuration (see \S\ref{S3.1}) on $(S^2)^n$.  Then
there is
a unique choice of initial point $x=x(r,\xi)$ such that $P$ closes up.
\end{thm}

\begin{rem}
Let $\t{\nu}$ be the atomic ``measure'' on $S^2$ which assigns mass $r_i$ to the
point $\xi_i$, $1 \leq i \leq n$, \emph{keeping track of the order of the
$\xi_i$'s}.  Then the rule that assigns $x=x(\t{\nu})=x(\xi,r)$ above is
$PSL_2(\C)$-equivariant and is a multiplicative analogue of the conformal
center of mass, $C(\nu)$, of Douady and Earle \cite{DE}, see also
\cite[\S4]{MZ}.
Here $\nu$ is the measure $\nu = \sum_{i=1}^n r_i \delta(\xi - \xi_i)$.
\end{rem}

\begin{rem}
We may use Theorem 1.3 to construct a length-shrinking flow on $CPol_n$.
Namely, let $0 \leq t \leq 1$.  Replace the weights $r=(r_1,...,r_n)$ by
$tr=(tr_1,...,tr_n)$.  We have
$$
M_r \xrightarrow{\ga_r} Q_{sst}(r) \cong Q_{sst}(tr) \xleftarrow{\ga_{tr}}
M_{tr}.
$$
\end{rem}
The composition $\ga_{tr}^{-1}\circ \ga_r$ is the length-shrinking flow.  Note
that $Q_{sst}(r)$ and $Q_{sst}(tr)$ are canonically isomorphic as complex
analytic spaces.  We obtain a curve $x(t)=x(tr, \xi)$.  We have

\begin{thm}
$\lim_{t \to 0} x(tr, \xi)=C(\nu)$, the conformal center of mass of Douady and
Earle.
\end{thm}

\begin{rem}
We see that $C(\nu)$ is ``semi-classical,'' it depends only on the limit as the
curvature goes to zero (or the speed of light goes to infinity), see
\S\ref{S3.3}.
\end{rem}

Our final theorems are connected with the study of certain integrable systems
on $M_r$ obtained by ``bending an n-gon along nonintersecting diagonals''
Precisely, we proceed as follows.  We define the diagonal $d_{ij}$ of $P$ to be
the geodesic segment joining $x_i$ to $x_j$.  Here we assume $i<j$.  We let
$\ell_{ij}$ be the length
of $d_{ij}$.  Then $\ell_{ij}$ is a continuous function on $M_r$ but is not
smooth at the points where $\ell_{ij}=0$.  We have the following description of
the Hamiltonian flow of $\ell_{ij}$ (it is defined provided $\ell_{ij} \ne 0$).

\begin{thm}
The Hamiltonian flow $\Psi_{ij}^t$ of $\ell_{ij}$ applied to an n-gon $P \in
M_r$ is obtained as follows.  The diagonal $d_{ij}$ separates $P$ into two
halves.  Leave one half fixed and rotate the other half at constant speed 1
around $d_{ij}$.
\end{thm}

\no For obvious reasons we call $\Psi_{ij}^t$ ``bending along $d_{ij}$.''

\begin{defn}
We say two diagonals $d_{ij}$ and $d_{ab}$ of $P$ do not intersect if the
interiors of $d^*_{ij}$ and $d^*_{ab}$ do not intersect, where $d^*_{ij}$ (resp.
$d^*_{ab}$) is the diagonal of a convex planar n-gon $P^*$ corresponding to
$d_{ij}$
(resp. $d_{ab}$).
\end{defn}
\no We then have

\begin{thm}
\label{Fla}
Suppose $d_{ij}$ and $d_{ab}$ do not intersect, then
$$
\{\ell_{ij},\ell_{ab}\}=0.
$$
\end{thm}

\vspace{-5mm}
\begin{rem}
We give two proofs of this theorem.  The first is a direct computation of the
Poisson bracket due to Hermann Flaschka.  The second is an elementary geometric
one depending on the description of the flows in Theorem 1.10.  It corresponds
to the geometric intuition that we may wiggle flaps of a folded piece of paper
independently if the fold lines do not intersect.
\end{rem}

We obtain a maximal collection of commuting flows if we draw
a maximal collection of nonintersecting diagonals $\{d_{ij}, (i,j) \in
I\}$. Later we will take the collection of all diagonals starting at the first vertex, $I=\{(1,3),(1,4),...,(1,n-1)\}$.  Each such collection
corresponds to a triangulation of a fixed convex planar n-gon $P^*$.  There are
$n-3$ diagonals in such a maximal collection.  Since dim $M_r =2n-6$, we obtain
\begin{thm}
For each triangulation of a convex planar n-gon $P^*$ we obtain an integrable
system on $M_r$.  Precisely, we obtain a Hamiltonian action of an (n-3)-torus
on $M_r$ which is defined on the Zariski open subset $M'_r$ defined by the
nonvanishing of the lengths of the diagonals in the triangulation.
\end{thm}

We have a simple description of the angle variables and the momentum
polyhedron attached to the above integrable system.  Let $M_r^o \subset M'_r$
be the subset such that none of the $n-2$ triangles in the triangulation are
degenerate.  Let $\widehat{\th}_{ij}$ be the dihedral angle at $d_{ij}$.  Put
$\th_{ij}=\pi-\widehat{\th}_{ij}$.  Then the $\th_{ij}$ are angle
variables.

To obtain the momentum polyhedron we follow \cite{HK} and note that there are
three triangle
inequalities associated to each of the $n-2$ triangles in the triangulation.
These are linear inequalities in the $\ell_{ij}$'s and the $r_{ij}$'s.  If they
are satisfied, we can build the $n-2$ triangles then glue them together and get
an n-gon $P$ with the required side-lengths $r_i$ and diagonal lengths
$\ell_{ij}$.  We obtain

\begin{thm}
The momentum polyhedron of the above torus action (the image of $M_r$ under the
$\ell_{ij}$'s) is the subset of $(\R_{\geq 0})^{n-3}$ defined by the $3(n-2)$
triangle inequalities above.
\end{thm}

\no As a consequence we obtain

\begin{cor}
The functions $\ell_{ij}, \, (i,j)\in I$, are functionally independent.
\end{cor}

Our results on $n$-gon linkages in $\h^3$ are the analogues of those of
\cite{KM2} for $n$-gon linkages in $\E^3$.  We conclude the paper by comparing
the symplectic manifolds $M_r(\h^3)$ and $M_r(\E^3)$.  Assume henceforth that
$r$ is not on a wall of $D_n$.

Since the Euclidean Gauss map $\ga_e: M_r(\h^3) \to Q_{sst}(r)$ is a canonical
diffeomorphism as is the hyperbolic Gauss map $\ga_h: M_r(\E^3) \to Q_{sst}(r)$
we obtain

\begin{thm}
The hyperbolic and Euclidean Gauss maps induce a canonical diffeomorphism
$$
M_r(\E^3) \simeq M_r(\h^3).
$$
\end{thm}

The last part of the paper is devoted to proving

\begin{thm} \label{1.17}
$M_r(\E^3)$ and $M_r(\h^3)$ are (noncanonically) symplectomorphic.
\end{thm}

This theorem is proved as follows.  Let $X_\ka$ be the complete simply-connected
Riemannian manifold of constant curvature $\ka$.  In \cite{Sa}, Sargent proved
that there exists $\al > 0$ and an analytically trivial fiber bundle $\pi : \e
\to (-\8,\al)$ such that $\pi^{-1}(\ka) = M_r(X_\ka)$.  We construct a closed
relative 2-form $\om_\ka$ on $\e|_{(-\8,0]}$ such that $\om_\ka$ induces a
symplectic form on each fiber of $\pi$ and such that the family of cohomology
classes $[\om_\ka]$ on $\e|_{(-\8,0]}$ is parallel for the Gauss-Manin
connection.  Theorem \ref{1.17} then follows from the Moser technique \cite{M}.

The results are closely related to but different from those of \cite{GW} and
\cite{A}.
\medskip

{\sc Acknowledgments.}
It is a pleasure to thank Hermann Flaschka for his help and
encouragement.  He
explained to us the set-up for the Sklyanin bracket (see \S\ref{S4.1}) and
provided us
with the first proof of Theorem \ref{Fla}.  Also, this paper was inspired by
reading
\cite{FR} when we realized that the dressing action of $SU(2)$ on $B^n$ was
just the natural action of $SU(2)$ on based hyperbolic n-gons.  We would also
like to thank Jiang-Hua Lu for explaining the formulas of \S\ref{form}
to us.  We would also like to thank her for pointing out that it was
proved in \cite{GW} that the cohomology class of the symplectic forms
$\om_\eps$ on an adjoint orbit in the Lie algebra of a compact group
was constant.

\section{Criteria for the moduli spaces to be smooth and nonempty} \label{S2}

In this chapter we will give necessary and sufficient conditions for the moduli
space $M_r$ to be nonempty and sufficient conditions for $M_r$ to be a smooth
manifold.

First we need some more notation.  Let $*$ be the point in $\h^3$ which is fixed
by $PSU(2)$.  We let $Pol_n(*)$ denote the space of n-gons $[x_1,...,x_{n+1}]$
with
$x_1=*$ and $CPol_n(*)$ = $CPol_n \cap Pol_n(*)$.  We let $\t{N}_r \subset
Pol_n(*)$ be the subspace of those n-gons $P=[x_1,...,x_{n+1}]$ such that
$d(x_i,x_{i+1})=r_i$, $1 \leq i \leq n$.  We put $N_r=\t{N}_r/PSU(2)$ and
$\t{M}_r = \t{N}_r \cap CPol_n(*)$.  Hence, $M_r = \t{M}_r/PSU(2)$.

Let $\pi:CPol_n \to (\R_{\geq 0})^n$ be the map that assigns to an n-gon $e$ its
set of side-lengths.  $\pi(e)=(r_1,...,r_n)$ with $r_i = d(x_i,x_{i+1}), \; 1
\leq i \leq n$.

\begin{lem}
The image of $\pi$ is the closed polyhedral cone $D_n$ defined by the
inequalities
$$
r_1 \geq 0, ... ,r_n \geq 0
$$
and the triangle inequalities
$$
r_i \leq r_1 + \cdots + \hat{r}_i + \cdots + r_n, \; 1 \leq i \leq n
$$
(here the \, $\hat{ }$ means that $r_i$ is omitted).
\end{lem}

\proof The proof is identical to the proof of the corresponding statement for
Euclidean space, \cite[Lemma 1]{KM1}.  \qed

\medskip
We next give sufficient conditions for $M_r$ to be a smooth manifold.  We will
use two results and the notation from $\S\ref{S4.3}$ (the reader will check that
no
circular reasoning
is involved here).  By Theorem \ref{sq} we find that $M_r$ is a symplectic
quotient.
$$
M_r \cong (\varphi|_{\t{N}_r})^{-1}(1)/SU(2)
$$
By Lemma \ref{reg}, 1 is a regular value of $\varphi$ unless there exists $P \in
\t{M}_r$ such that the infinitesimal isotropy $(su_2)|_P = \{x\in SU(2):
\hat{X}(P)=0 \}$ is nonzero.

\begin{defn}
An n-gon P is degenerate if it is contained in a geodesic.
\end{defn}

We now have
\begin{lem}
$M_r$ is singular only if there exists a partition $\{1,...,n\} = I \amalg J$
with $\#(I) > 1, \#(J) > 1$ such that
$$
\sum_{i \in I} r_i = \sum_{j \in J} r_J.
$$
\end{lem}

\proof Clearly $(su_2)|_P=0$ unless $P$ is degenerate.  But if $P$ is
degenerate there exists a partition $\{1,...,n\} = I \amalg J$ as above ($I$
corresponds to the back-tracks and $J$ to the forward-tracks of $P$). \qed

\begin{rem}
In the terminology of \cite{KM1}, \cite{KM2}, $M_r$ is
smooth unless $r$ is on a wall of $D_n$.  Note that if $|I|=1$ or $|J|=1$ then
$r \in \partial D_n$ and $M_r$ is reduced to a single point.
\end{rem}

There is a technical point concerning smoothness.  We could also define $M_r$
as the fiber of $\bar{\pi}: CPol_n/PSL_2(\C) \to D_n$ over $r$.  It is not
quite immediate that smoothness of the symplectic quotient coincides with the
smoothness of $\bar{\pi}^{-1}(r).$  Fortunately, this is the case (note $r$ is a
regular value of $\bar{\pi} \: \Leftrightarrow \: r$ is a regular value of
$\pi)$.

\begin{lem} \label{regval}
$r$ is a regular value of $\pi \: \Leftrightarrow \: 1$ is a regular value of
$\varphi| \t{N}_r$.
\end{lem}

\proof The lemma follows from a consideration of the diagram
$$
\xymatrix{
& \t{N}_r \ar[d] \ar[dr]^{\varphi_r} \\
CPol_n(*) \ar[r] \ar[dr]_{\pi} & Pol_n(*) \ar[r]^{\varphi} \ar[d] & B \\
& (\R_{\geq 0})^n  }
$$
and the observation that $\varphi:Pol_n(*) \to B$ (see \S\ref{S4.2}) and the
side-length map
$Pol_n(*) \to (\R_{\geq 0})^n$ are obviously submersions.  Here we have
abbreviated $\varphi|\t{N}_r$ to $\varphi_r$. \qed

\section{The geometric invariant theory of hyperbolic polygons} \label{S3}

\subsection{The hyperbolic Gauss map and weighted quotients of the configuration
spaces of points on the sphere} \label{S3.1}

The goal of the next two sections is to construct a natural
homeomorphism $\gamma : M_r \to Q_{sst}$ where $Q_{sst}$ is the
$r$-th weighted quotient of $(S^2)^n$ by $PSL_2(\C)$ constructed
in \cite{DM} in the case that $M_r$ is smooth.  $Q_{sst}$ is a
complex analytic space.  We now review the construction of
$Q_{sst}$.

Let $M \subset (S^2)^n$ be the set of n-tuples of distinct points. Then
$Q=M/PSL_2(\C)$ is a (noncompact) Hausdorff manifold.

\begin{defn}
A point $\vec{u} \in (S^2)^n$ is called r-stable (resp. semi-stable) if
$$
\sum_{u_j=v} r_j < \frac{|r|}{2} \:(\textnormal{resp.} \leq \frac{|r|}{2})
$$
for all $v \in S^2$. Here $|r| = \sum_{j=1}^n r_j$.  The set of stable and
semi-stable points will be
denoted by $M_{st}$ and $M_{sst}$ respectively.  A semi-stable point $\vec{u}
\in (S^2)^n$ is said to be a nice semi-stable point if it is either stable
or the orbit $PSL_2(\C) \vec{u}$ is closed in $M_{sst}$.
\end{defn}
We denote the space of nice semi-stable points by $M_{nsst}$.  We have the
inclusions
$$
M_{st} \subset M_{nsst} \subset M_{sst}.
$$
Let $M_{cusp} = M_{sst} - M_{st}$. We obtain the points $M_{cusp}$ in the
following way.  Partition $S=\{1,...,n\}$ into disjoint sets $S=S_1 \cup S_2$
with $S_1=\{i_1,...,i_k\}, \: S_2=\{j_1,...,j_{n-k}\}$ in such a way that
$r_{i_1}
+ \cdots + r_{i_k} = \frac{|r|}{2}$ (whence $r_{j_1} + \cdots + r_{j_{n-k}} =
\frac{|r|}{2}$).  Then
$\vec{u}$ is in $M_{cusp}$ if either $u_{i_1}= \cdots = u_{i_k}$ or $u_{j_1}=
\cdots
= u_{j_{n-k}}$.  The reader will verify that $\vec{u} \in M_{cusp}$ is a nice
semi-stable point if and only if both sets of the equations above hold.
relation $\r$ via: \\
$\vec{u} \equiv \vec{w} \pmod{\r}$ if either \\
(a) $\vec{u},\vec{w} \in M_{st}$ and $\vec{w} \in PSL_2(\C)\vec{u}$, \\
or \\
(b) $\vec{u},\vec{w} \in M_{cusp}$ and the partitions of $S$ corresponding to
$\vec{u},\vec{w}$ coincide.

The reader will verify that $\vec{u},\vec{w} \in M_{nsst}-M_{st}$ then $\vec{u}
\equiv
\vec{w}
\pmod{\r}$ if and only if $\vec{w} \in PSL_2(\C)\vec{u}$.

It is clear that $\r$ is an equivalence relation. Set
$$
Q_{sst} = M_{sst}/\r, \;Q_{nsst} = M_{nsst}/\r, \;Q_{st} =
M_{st}/\r, \;Q_{cusp} = M_{cusp}/\r
$$
each with the quotient topology.  The elements of $Q_{cusp}$ are uniquely
determined by their partitions.  Thus $Q_{cusp}$ is a finite set.  It is clear
that each equivalence class in $Q_{cusp}$ contains a unique $PSL_2(\C)$-orbit of
nice semi-stable points whence the inclusion
$$
M_{nsst} \subset M_{sst}
$$
induces an isomorphism
$$
Q_{nsst} = M_{nsst}/PSL_2(\C) \to Q_{sst}.
$$

In case $r_1,...,r_n$ are rational then the quotient space
$Q_{sst}$ can be given a structure of a complex projective variety
by the techniques of geometric invariant theory applied to certain
equivariant projective embeddings of $(S^2)^n$, see
\cite[\S4.6]{DM}.  This concludes our review of \cite[\S4]{DM}. We
now establish the connection with the moduli space $M_r$.

For the rest of this section (with one exception) we will use the
ball model of hyperbolic space (so $*=(0,0,0)$).  We will
compactify $\h^3$ by enlarging the open three ball to the closed
three ball, thus we add $S^2=\partial_\8 \h^3$. Each point of
$S^2$ corresponds to an equivalence class of geodesic rays in
$\h^3$.  Two rays $\al$ and $\be$ are equivalent if they are
asymptotic, i.e. $\lim_{t \to \8} \al (t) = \lim_{t \to \8} \be
(t)$ in the closed three ball. Intrinsically the equivalent rays
are characterized by the property that they are within finite
Hausdorff distance from each other.

In what follows all geodesic segments, geodesics and geodesic rays
will be parameterized by arc-length. We now define the hyperbolic
Gauss map $\ga$ (in various incarnations). Let $\si=[x,y],
\;x,y\in \h^3$, be the oriented geodesic segment from $x$ to $y$.
Let $\t{\si}(0)$ be the ray, $\t{\si}:[0,\8) \to \h^3$ with
$\t{\si}(0)=x$ and $\t{\si}(\ell)=y$ (here $\ell = \ell(\si)$ is
the length of the geodesic segment $\si$).  We define the
(forward) Gauss map $\ga$ on oriented segments by $$ \ga(\si) =
\lim_{t\to\8} \t{\si}(t). $$ We may now define $\ga: \t{N}_r \to
(S^2)^n$ by $$ \ga(e)=(\ga(e_1),...,\ga(e_n)). $$ One of the main
results of this paper is the following theorem -- an analogue for
Poisson actions of the theorem of Kirwan, Kempf, and Ness,
\cite{Ki}, \cite{KN}.

\begin{thm}
\label{KKN}
\begin{itemize}
\item[(i)] $\ga(\t{M}_r)\subseteq M_{nsst}$.
\item[(ii)] If $P$ is nondegenerate, then $\ga(P) \in M_{st}$.
\item[(iii)] $\ga$ induces a real analytic homeomorphism $\ga:M_r \to
Q_{sst}$.
\item[(iv)] $M_r$ is smooth if and only if $M_{st} = M_{sst}$.  In this case
$Q_{sst}$ is also smooth and $\ga: M_r \to Q_{sst}$ is an analytic diffeomorphism.
\end{itemize}
\end{thm}

Let $\eta\in S^2$.  We recall the
definition of the geodesic flow $\phi_{\eta}^t$ associated to $\eta$.
(Strictly speaking,
this flow is rather the projection to $\h^3$ of the restriction of the geodesic flow on $UT(\h^3)$
to the stable submanifold corresponding to $\eta$.) Given
$z \in \h^3$ there is a unique arc-length parameterized ray $\si$ with
$\si(0)=z, \; \lim_{t\to\8} \si(t)=\eta$.  By definition,
$$
\phi_{\eta}^t(z)=\si(t).
$$

We will also need the definition of the Busemann function $b(x,\xi), \; x\in
\h^3, \xi \in \partial_\8 \h^3$.  Let $\si$ be an arc-length parameterized
geodesic ray from $*$ to $\xi$.  Then
$$
b(x,\xi)= \lim_{t\to\8}(d(x,\si(t))-t).
$$
Note that for $k \in PSU(2) = Stab(*)$ we have
$$
b(kx, k\xi) = b(x,\xi).
$$
Also, in the upper half space for $\h^3$, we have
$$
b((x,y,z),\8)=-\log z.
$$

\no We now prove

\begin{lem}
For fixed $\xi$, $-\nabla b(x,\xi)$ is the infinitesimal generator of the
geodesic flow $\phi_\xi^t$.
\end{lem}
\proof  From the first formula above, it suffices to check this statement in
the upper half space model (so $*=(0,0,1)$) for $\xi = \8$.  By the second
formula
$$
-\nabla b(x,y,z)= z \frac{\partial}{\partial z}\ .
$$
\qed

\no We will now prove (i) and (ii) in the statement of Theorem \ref{KKN}
above.

\begin{lem}
\begin{itemize}
\item[(i)] $\ga(\t{M}_r) \subseteq M_{nsst}$.
\item[(ii)] $P$ is nondegenerate $\Leftrightarrow \; \ga(P) \in M_{st}.$
\end{itemize}
\end{lem}
\proof Let $P\in\t{M}_r$ be a polygon with the vertices $x_1,..., x_{n+1}= x_n$, we will use the notation
$x_i(t), 0\leq t\leq r_i,$ for the parameterized edge $e_i$ (so that $x_i(0)=x_i$).
We test stability of $\ga(P)$ with respect to a point $\eta\in S^2$. Let
$b(x):= b(x,\eta)$ be the corresponding Busemann function. Then for any unit vector
$v\in T_x(\h^3)$
\begin{equation}
\label{g1}
-\nabla b(x)\cdot v \leq 1
\end{equation}
with the equality if and only if the geodesic ray $\exp(\R_+ v)$ is asymptotic to $\eta$.
Similarly,
\begin{equation}
\label{g2}
-\nabla b(x)\cdot v \ge -1
\end{equation}
with the equality if and only if the geodesic ray $\exp(\R_- v)$ is asymptotic to $\eta$.
Let $I \subset \{1,...,n\}$ be the subset of indices such that $\ga(e_i)=\eta$.
Let $J$ be the complement of $I$ in $\{1,...,n\}$.
Put $r_I = \sum_{i\in I} r_i, \; r_J = \sum_{j\in J} r_j$.
Since the polygon $P$ is closed, using (\ref{g1}) and (\ref{g2})  we get:
\begin{eqnarray*}
0= -b \bigg|_{x_1}^{x_{n+1}}= \sum_{i=1}^{n-1} \int_{0}^{r_i}- x_i'(t) \cdot \nabla b(x_i(t))dt \\
\geq \sum_{i\in I} r_i - \sum_{j\in J} r_j= r_I - r_J
\end{eqnarray*}
with the equality if and only if every edge $e_j, j\in J,$ is contained in
the geodesic through $\eta$ and $x_j$.
Thus $r_I\leq r_J$, i.e. $\ga(P)$ is semi-stable. If
$\ga(P)$ is not stable then each
edge $e_i, 1\leq i\leq n,$ of $P$ is contained in the
geodesic through $\eta$ and $x_i$, which implies
that this geodesic is the same for all $i$. Hence $P$
is degenerate in this case.  \qed

\medskip
In order to prove that $\ga: M_r\to Q_{nsst}$ is injective
and surjective, we will first
need to study a certain dynamical system $f_{r,\xi} \in$Diff$(\h^3)$
attached to the configuration of $n$ points $\xi = (\xi_1,...,\xi_n)$
on $S^2$ weighted by $r=(r_1,...,r_n)$.
The weights $r$ will usually be fixed and we will drop $r$ in $f_{r,\xi}$.

\subsection{A dynamical system on $\h^3$ and the proof that the Gauss map is an
isomorphism}

Let $\xi=(\xi_1,...,\xi_n) \in (S^2)^n$.  We define a diffeomorphism
$f_\xi:\h^3 \to \h^3$ as follows.  Assume that $r=(r_1,...,r_n)\in D_n$ is
given.  Let $z\in \h^3$ be given.  Let $\si_1$ be the ray emanating from $z$
with $\lim_{t\to\8} \si_1(t)= \xi_1$.  Put $x_1=z$ and $x_2 = \si_1(r_1)$.
Now let $\si_2$ be the ray emanating from $x_2$
with $\lim_{t\to\8} \si_2(t)= \xi_2$.  Put $x_3 = \si_2(r_2)$.
We continue in this way until we obtain $x_{n+1}=\si_n(r_n)$ where
$\si_n$ is the geodesic ray emanating from $x_n$ with $\lim_{t\to\8} \si_n(t)=
\xi_n$.  We define $f_\xi$ by $f_\xi(z) =
x_{n+1}$. Note that the polygon $P=(x_1,...,x_{n+1})$ belongs to $\tilde{N}_r$.

We now give another description of $f_\xi$:
$$
f_\xi = \phi_{\xi_n}^{r_n} \circ \cdots \circ \phi_{\xi_1}^{r_1}
$$
where $\phi_\xi^t$ is the time $t$ geodesic flow towards $\xi$.
We may interpret the previous formula for $f_\xi$ as a product (or
multiplicative) integral \cite{DF}.  Partition the interval $[0,1]$ into
$n$ equal subintervals, $0=t_0<t_1< \cdots <t_n=1$.  Let $\nu$ be the atomic
measure on $[0,1]$ given by $\nu(t)=\sum_{i=0}^{n-1} r_{i+1} \del(t-t_i)$.
Let $\la:[0,1] \to S^2$ be the map given by $\la|[t_i, t_{i+1}) = v_i, \; 0
\leq i \leq n-1$.  Define $A:[0,1]\to C^\8(\h^3, T(\h^3))$ by $A(t)(z) =
\nabla \,b(z,\la(t))$.  Then in the
notation of \cite{DF},
$$
f_\nu = \prod_0^{1} e^{A(t)d\nu(t)}.
$$
in Diff$(\h^3)$.

\begin{rem} In fact, in \cite{DF} the only integrals considered take values in
$G\hspace{-1pt}L_n(\C)$.  
We have included the above formula to stress the analogy with the
conformal center of mass.  The above integral is the multiplicative analogue of
the gradient of the averaged Busemann function
$$
\nabla b_\nu(z) = \int_{S^2} \nabla b(z,\eta)(\la_* d\nu)(\eta)
$$
used to define the conformal center of mass (see \S\ref{S3.3}).
\end{rem}

\no We will first prove

\begin{prop}
Suppose $\xi$ consists of three or more distinct points.  Then $f_\xi$ is a
strict contraction.  Hence, if $\xi$ is stable, $f_\xi$ is a strict
contraction.
\end{prop}

\no We will need the following lemma

\begin{lem}
Let $\xi \in S^2$ and $\phi_{\xi}^t$ be the geodesic flow towards $\xi$.  Then,
for each $t>0$,
\begin{itemize}
\item[(i)] $d(\phi_{\xi}^t(z_1),\phi_{\xi}^t(z_2)) \leq d(z_1,z_2)$ with
equality
if and only if $z_1$ and $z_2$ belong to the same geodesic $\eta$ with
end-point $\xi$.

\item[(ii)] If $Z \in T_z(\h^3)$ is a tangent vector, then
$$
||D \phi_{\xi}^t(Z)|| \leq ||Z||
$$
with equality if and only if $Z$ is tangent to the geodesic $\eta$ through $z$
which is asymptotic to $\xi$.
\end{itemize}
\end{lem}
\proof  We prove (ii) noting that (ii) implies (i).  Use the upper half
space model for
$\h^3$ and send $\xi$ to $\8$.  Then, if $Z=(a,b,c)$ is tangent to $\h^3$ at
$(x,y,z)$,
we have
$$
\phi_{\8}^t(x,y,z) = (x,y,e^t z)
$$
and
$$
||D \phi_{\8}^t(a,b,c)|| \bigg|_{(x,y,e^t z)} = \sqrt{e^{-2t}a^2 +e^{-2t}b^2 +
c^2}\ .
$$
\qed

\begin{rem}
1. In case of equality in (i), the points
$z_1,\,z_2, \, \phi_{\xi}^t(z_1), \, \phi_{\xi}^t(z_2)$ all
belong to $\eta$.

2. The above lemma also follows from the fact that
$D \phi_{\xi}^t(Z)$ is a stable Jacobi field along $\eta$.
\end{rem}

We can now prove the proposition.  By the previous lemma, $f_\xi =
\phi_{\xi_n}^{r_n} \circ \cdots \circ \phi_{\xi_1}^{r_1}$ does not increase
distance.  Suppose then that $d(f_{\xi}(z_1),f_{\xi}(z_2)) =
d(z_1,z_2)$.  Then,
$$
d(\phi_{\xi_1}^{r_1} (z_1),\phi_{\xi_1}^{r_1} (z_2)) = d(z_1,z_2).
$$
Hence, $z_1,\; z_2, \; \phi_{\xi_1}^{r_1} (z_1), \; \phi_{\xi_1}^{r_1} (z_2)$
are
all on the geodesic $\eta_1$ joining $z_1$ to $\xi_1$.

Next,
$$
d(\phi_{\xi_2}^{r_2}(\phi_{\xi_1}^{r_1}(z_1)),\phi_{\xi_2}^{r_2}
(\phi_{\xi_1}^{r_1} (z_2))) =
d(\phi_{\xi_1}^{r_1}(z_1),\phi_{\xi_1}^{r_1}(z_2)).
$$
Hence, $\phi_{\xi_1}^{r_1}(z_1), \; \phi_{\xi_1}^{r_1}(z_2), \;
\phi_{\xi_2}^{r_2}(\phi_{\xi_1}^{r_1}(z_1)),\;\phi_{\xi_2}^{r_2}(\phi_{\xi_1}^{r
_1}(z_2))$
are all on the same geodesic.  This geodesic is necessarily $\eta_2$, the
geodesic joining $\phi_{\xi_1}^{r_1}(z_1)$ to $\xi_2$, since it contains
$\phi_{\xi_1}^{r_1}(z_1)$ and $\phi_{\xi_2}^{r_2}(\phi_{\xi_1}^{r_1}(z_1))$.
But since $\eta_2$ contains $\phi_{\xi_1}^{r_1}(z_1)$ and
$\phi_{\xi_1}^{r_1}(z_2)$ it also coincides with $\eta_1$.  Hence, either
$\xi_2=\xi_1$ or $\xi_1$ is the opposite end $\check{\xi}_1$ of the geodesic
$\eta_1$.  We continue in this way and find that either $\xi_i=\xi_1$ or $\xi_i
= \check{\xi}_1$, for all $1 \leq i \leq n$. \qed

\medskip
Our next goal is to prove that $f_\xi$ has a fixed-point in $\h^3$.  Let $H$ be
the convex hull of $\{\xi_1,...,\xi_n\}$.  Let $\be_i$ be the negative of the
Busemann function associated to $\xi_i$ so $\be_i$ \emph{increases} along
geodesic rays directed toward $\xi_i$.
Fix a vector $r=(r_1,...,r_n)$ and $r$-stable configuration
$\xi=(\xi_1,...,\xi_n)
\in S^2$.  For $1 \geq h > 0$, we shall consider $f_{hr,\xi}:\h^3 \to \h^3$
where $hr = (hr_1,...,hr_n)$. Note that $f_{hr,\xi}(H)\subset H$.

\begin{lem}
\label{Bus}
There exist open horoballs, $\O_i, \; 1\leq i \leq n$, centered at $\xi_i$,
which depend only on $r$ and $\xi$, such that for each $1 \leq i \leq n$, if
$x \in \O_i \cap H$, then
$$
\be_i(f_{hr,\xi}(x))<\be_i(x)
$$
(so $f_{hr,\xi}(x)$ is ``further away from'' $\xi_i$ than $x$).
\end{lem}

\proof The angle between any two geodesics asymptotic to $\xi_i$
is zero, thus by continuity, for each $\eps>0$, there exists a
horoball $\O_i(\eps)$ centered at $\xi_i$ so that for each $x\in
\O_i(\eps)\cap H$ and for each point $\xi_j$ which is different
from $\xi_i$, the angle between the geodesic ray from $x$ to
$\xi_j$ and $\nabla b(x,\xi_i)$ is $\leq \eps$. Let $I=\{\ell\in
\{1,...,n\}:\xi_\ell=\xi_i\}, \;J:=\{1,...,n\}-I$.  Recall the
stability condition means: $$ r_I:=\sum_{\ell \in I} r_\ell < r_J
:= \sum_{j\in J}r_j $$ thus we can choose $\pi/2 > \eps>0$ so that
$$ r_I - \cos(\eps)r_J<0. $$ We define $L$ by $\O_i(\eps)= \{\be_i
> L\}$ for this choice of $\eps$, then we define $\O_i$ by $\O_i:=
\{ \be_i > L + |r|\}$. Pick $x \in H \cap \O_i$, this point is the
initial vertex of the linkage $P$ with vertices $$ x_1=x,\;
x_2=\phi_{\xi_1}^{hr_1}(x_1),...,\;
x_{n+1}=\phi_{\xi_n}^{hr_n}(x_n)=f_{hr,\xi}(x). $$ Note that since
the length of $P$ equals $h|r|$, (and $h \leq 1$), the whole
polygon $P$ is contained in $H \cap \O_i(\eps)$.  We let $x_j(t),
t\in [0, hr_j]$ be the geodesic segment connecting $x_j$ to
$x_{j+1}$ (parameterized by the arc-length). Then,
\begin{eqnarray*}
\be_i(x_{n+1})- \be_i(x_1) & = & \sum_{k=1}^n \be_i(x) \bigg|_{x_k}^{x_{k+1}}\\
& = & \sum_{k=1}^n \int_0^{hr_k} \nabla \be_i(x_k(t)) \cdot x'_k(t)dt\\
& = & \sum_{\ell \in I} \int_0^{hr_\ell} \nabla \be_i(x_\ell(t)) \cdot
x'_\ell(t)dt
+ \sum_{j \in J} \int_0^{hr_j} \nabla \be_i(x_j(t)) \cdot x'_j(t)dt.
\end{eqnarray*}
Recall that $||\nabla \be_i(x_k(t))||=1, \; ||x'_k(t)||=1$,
if $\ell\in I$ then
$$
\nabla \be_i(x_\ell(t))\cdot x'_\ell(t)=1,
$$
if $j \in J$ then
$$
\nabla \be_i(x_j(t))\cdot x'_j(t) \leq -\cos(\eps)
$$
since $x_j(t) \in \O_i(\eps)$ for each $ 0 \leq t \leq hr_j$.  Thus,
$$
\be_i(x_{n+1}) - \be_i(x_1) \leq h \sum_{\ell \in I}r_{\ell} \ \ - h \cos(\eps) \sum_{j\in
J}r_j = h(r_I-\cos(\eps)r_J)<0.
$$
\qed

\no We let $\O'_i:= \{\be_i > L + 2|r|\}$, then

\begin{prop}
$f_{hr,\xi}$ has a fixed point in $K:=H -\bigcup_{i=1}^n \O'_i$.
\end{prop}
\proof  We claim that if $x \in H-\bigcup_{i=1}^n \O_i$ then for all $m \geq
0, \; f_\xi^{(m)}(x) \notin \bigcup_{i=1}^n \O'_i$.  We first treat the case
$m=1$.  Since $d(x, f_\xi(x)) \leq |r|$ we see that $x \in H-\bigcup_{i=1}^n
\O_i$ implies $f_\xi(x) \notin \bigcup_{i=1}^n \O'_i$.  But if there exist an
$m-1$ such that $y = f_\xi^{(m-1)}(x) \in \bigcup_{i=1}^n (\O_i-\O'_i)$, then
$f_\xi^{(m)}(x)=f_\xi(y) \notin \bigcup_{i=1}^n \O'_i$ by Lemma \ref{Bus}
and the claim is proved.

We find that the sequence $\{f_{\xi}^{(m)}(x)\}$ is relatively compact and contained
in $K$. Let $A\subset K$ be the accumulation set for this sequence. This is a compact
subset such that $f_{\xi}(A)\subset A$. If $f_\xi$ does not have a fixed point in $A$
then the continuous function $\theta(x):= d(x, f_{\xi}(x)), x\in A$ is bounded away from zero.
Let $x_0\in A$ be a point where $\theta$ attains its minimum. However (since $f_\xi$ is a strict contraction)
$$
\theta(f_\xi (x_0))= d(f_\xi(x_0), f^2_\xi(x_0))< d(x_0, f(x_0))= \theta(x_0),
$$
contradiction. \qed

We can now prove Theorem \ref{KKN}.
We first prove that $\ga: M_r\to Q_{sst}$ is
injective.  This easily reduces to proving that if $P,Q \in \t{M}_r$ with
$\ga(P)=\ga(Q)$, then $P=Q$.  Let $x_1$ be the first vertex of $P$, $x'_1$ be
the first vertex of $Q$, and $\xi = \ga(P) = \ga(Q)$.  Since $P$ closes up, we
have $f_\xi(x_1)=x_1$.  Since $Q$ closes up, we have $f_\xi(x'_1)=x'_1$.  But,
$f_\xi$ is a strict contraction, hence $x_1=x'_1$.  It follows immediately that
$P=Q$.

We now prove that $\ga$ is surjective.  Let $\xi \in M_{st}$.  There
exists $x\in\h^3$ with $f_\xi(x)=x$.  Let $P$ be the n-gon with $\ga(P)=\xi$
and first vertex $x$.  Then $P$ closes up and we have proved that $\ga$ is onto
the stable points.  If $\xi$ is nice semi-stable but not stable, then
$\xi=\ga(P)$ for a suitable degenerate n-gon.  Hence, $\ga$ is surjective
and Theorem \ref{KKN} is proved.

\begin{rem}
We have left the proof that the inverse map to $\ga: M_r\to Q_{st}$
is smooth (resp. analytic) in the
case $M_r$ is smooth to the reader.  This amounts to checking that the
fixed-point of $f_\xi$ depends smoothly (resp. analytically) on $\xi$.
\end{rem}

\subsection{Connection with the conformal center of mass of Douady and Earle}
\label{S3.3}

In this section, we prove Theorem 1.8 of the Introduction.  We begin by
reviewing the definition of the conformal center of mass $C(\nu)\in \h^3$,
where $\nu$ is a stable measure on $S^2 = \partial_\8 \h^3$.  Here we are using

\begin{defn}
A measure $\nu$ on $S^2$ is stable if
$$
\nu(\{x\})< \frac{|\nu|}{2}, \; x\in S^2.
$$
Here, $|\nu|$ is the total mass of $\nu$.
\end{defn}

We define the averaged Busemann function, $b_\nu :\h^3\to \R$, by
$$
b_\nu(x)= \int_{S^2} b(x,\xi) d\nu(\xi).
$$

\no We recall the following proposition (\cite{DE}, \cite[Lemma 4.11]{MZ}):
\begin{prop}
Suppose $\nu$ is stable.  The $b_\nu$ is strictly convex and has a unique
critical point (necessarily a minimum).
\end{prop}

\begin{defn}
The conformal center of mass $C(\nu)$ is defined to be the above critical
point.  Thus,
$$
\nabla b_\nu |_{C(\nu)} = 0.
$$
\end{defn}

\no The main point is the following,
\begin{lem}
The assignment $\nu \to C(\nu)$ is $PSL_2(\C)$-equivariant,
$$
C(g_* \nu) = g C(\nu).
$$
Here $g_* \nu$ is the push-forward of $\nu$ by $g\in PSL(2,\C)$.
\end{lem}

We now return to the set-up of the previous sections.  We are given
$r=(r_1,...,r_n)$ and a stable configuration $\xi=(\xi_1,...,\xi_n) \in
(S^2)^n$.  We have the dynamical system $f_{tr,\xi}$ of the previous chapter,
with fixed-point $x = x(tr,\xi)$.  We put $\nu = \sum_{i=1}^n r_i \del(\xi -
\xi_i)$, where $\del$ is the Dirac probability measure supported on the origin in $\R^3$.
 We now have,
\begin{lem}
\label{5.6}
$$
\frac{d}{dt} f_{tr,\xi} \bigm|_{t=0} = -\nabla b_\nu\ .
$$
\end{lem}
\proof  We abbreviate $-\nabla b(x,\xi_i)$, the infinitesimal generator of the
geodesic flow associated to $\xi_i$, to $X_i$.  Thus we want to prove
$$
\frac{d}{dt} f_{tr,\xi}\bigm|_{t=0} = \sum_{i=1}^n r_i X_i.
$$
But if $\varphi_t$ and $\psi_t$ are flows with infinitesimal generators $X$ and
$Y$ respectively, then
\begin{eqnarray*}
  \frac{d}{dt}\biggm|_{t=0} \varphi \circ \psi(x) & = &
\frac{\partial^2}{\partial
  t_1 \partial t_2}\biggm|_{t_1=t_2=0} \varphi_{t_1} \circ \psi_{t_2}(x) \\
  & = & \frac{\partial}{\partial t_2}\biggm|_{t_2=0} \varphi_o \circ
\psi_{t_2}(x) +
  \frac{\partial}{\partial t_1}\biggm|_{t_1=0} \varphi_{t_1} \circ \psi_o(x) \\
  & = & \frac{\partial}{\partial t_2}\biggm|_{t_2=0} \psi_{t_2}(x) +
  \frac{\partial}{\partial t_1}\biggm|_{t_1=0} \varphi_{t_1}(x) \\
  & = & Y(x) + X(x).
\end{eqnarray*}
Recall that
$$
f_{tr,\xi}(x)= \phi^{tr_n}_{\xi_n} \circ \phi^{tr_{n-1}}_{\xi_{n-1}}
\circ \cdots \circ \phi^{tr_1}_{\xi_1}(x).
$$
Hence,
$$
\frac{d}{dt}\biggm|_{t=0} f_{tr,\xi}(x) = \sum_{i=1}^n X_i(x).
$$
\qed

We are ready to prove Theorem 1.8 of the Introduction.
We abbreviate $-\nabla b_\nu$ by $X$.
\begin{thm}
Let $x(tr,\xi)$ be the unique fixed point of $f_{tr,\xi}, \; 0 < t \leq 1$.
Then
$$
\lim_{t\to 0} x(tr,\xi) = C(\nu).
$$
\end{thm}
\proof  We note that $f_{0r,\xi}$ = id.  Hence, applying Lemma \ref{5.6}, the
Taylor approximation of $f_{tr,\xi}(x)$ around $t=0$ is
$$
f_{tr,\xi}(x) = x + tX(x) + t^2 R(x,t)
$$
(where $R(x,t)$ is smooth).
Let $\varphi (x,t)= X(x) + tR(x,t)$.  By definition, the conformal center of
mass $C(\nu)$ is the unique solution of
$$
\varphi(x,0)= X(x) = 0.
$$
Since $b_\nu$ is strictly convex, if $\nu$ is stable, \cite[Corollary 4.6]{MZ},
$C(\nu)$ is a nondegenerate zero of $X$ and we may apply the implicit function
theorem to solve
$$
\varphi(x,t)=0
$$
for $x$ as a function of $t$ near $(C(\nu),0)$.  Thus, there exists $\delta >0$
and smooth curve, $\hat{x}(t)$, defined for $|t|<\delta$, satisfying
\begin{itemize}
\item[(i)] $\varphi(\hat{x}(t), t)=0$
\item[(ii)] $\hat{x}(0)= C(\nu)$.
\end{itemize}
But clearly, (i) implies $f_{tr,\xi}(\hat{x}(t))=\hat{x}(t), \; 0 < t <
\delta$. Since the fixed-point of $f_{tr,\xi}, \; 0<t\leq1$, is unique, we
conclude $\hat{x}(t)=x(tr,\xi), \; 0<t<\delta$.  Hence,
$$
\lim_{t\to 0} x(tr,\xi)= \lim_{t \to 0} \hat{x}(t)=C(\nu).
$$
\qed

\section{The Symplectic geometry of $M_r(\h^3)$} \label{S4}

\subsection{The Poisson-Lie group structure on $B^n$} \label{S4.1}

In this section we let $G$ be any (linear) complex simple group,
$B=AN$ be the subgroup of the Borel subgroup such that $N$ is its unipotent
radical and $A$ is the connected component of the identity in a maximal split
torus over $\R$, and $K$ be a maximal compact subgroup. We
will construct a Poisson Lie group structure on $G$ which will
restrict to a Poisson Lie group structure on $B$. For the basic
notions of Poisson Lie group, Poisson action, etc. we refer the
reader to \cite{Lu}, \cite{GW}, \cite{LW}, and \cite{CP}.

Let $R_g$ and $L_g$ be the action of $g$ on $G$ by the right and left
multiplication respectively.  Let $\g$ denote the Lie algebra of $G$,  $\k$ be
the Lie algebra of $K$ and $\b$ be the Lie algebra of $B$. Then $\g= \b \oplus
\k$ and $G=BK$. Let $\rho_\k$ (resp. $\rho_\b$) be the projection on $\k$ (resp.
on $\b$). We define ${\r}:= \rho_\k - \rho_\b$ and let
$ \half \< , \>$ be the imaginary part of the Killing form on $\g$.
In the direct sum splitting $\g = \k  \oplus \b$ we see that $\k$ and $\b$ are
totally-isotropic subspaces dually paired by $\< , \>$.

Let $\varphi\in C^\8 (G)$. Define $D\varphi : G\to \g$ and
$D'\varphi : G\to \g$ by
$$
\< D'\varphi (g), \nu\> = \frac{d}{dt}|_{t=0} \varphi(g e^{t\nu})
$$
$$
\< D\varphi (g), \nu\> = \frac{d}{dt}|_{t=0} \varphi(e^{t\nu}g)
$$
for $\nu \in \g$.

We extend $\< , \>$ to a biinvariant element of  $C^\8  (G, S^2 T^*(G))$ again
denoted $\< , \>$. Now define $\nab \varphi \in C^\8  (G,  T^*(G))$ by
$$
\< \nab \varphi (g), x\> = d\varphi_g(x), x\in T_g(G)
$$
We have
$$
D\varphi (g) = dR^{-1}_g \nab \varphi(g)
$$
$$
D'\varphi (g) = dL^{-1}_g \nab \varphi(g)= Ad_{g^{-1}} D\varphi(g)
$$
The Sklyanin bracket $\{ \varphi, \psi\}$ is defined for $\varphi, \psi \in C^\8
(G)$ by
$$
\{ \varphi, \psi\}(g)= \half ( \<  {\r} D' \varphi(g), D' \psi(g)\> -
\< {\r}D \varphi(g), D\psi(g)\> )
$$
where ${\r}= \rho_\k -\rho_\b$.

\begin{lem}
The bracket $\{ \varphi, \psi\}$ is a Poisson bracket on $C^\8 (G)$.
\end{lem}
\proof See Theorem 1 of \cite{Se}.

\medskip
Let $w\in C^\8 (G, \La^2 T(G))$ be the bivector field corresponding to
$\{ \cdot , \cdot \}$. We now show that $\{ \cdot , \cdot \}$ induces a Poisson
bracket on $C^\8 (B)$.

\begin{lem}
$w(b)$ is tangent to $B$ , i.e. $w(b)\in \La^2 T_b (B) \subset \La^2 T(G)$ for
all $b\in B$.
\end{lem}
\proof It suffices to prove that if $\varphi$ vanishes identically on $B$ then
$\{ \varphi, \psi\}$ vanishes identically on $B$ for each $\psi$. However
if $\varphi$ vanishes identically on $B$ then $\nab \varphi(b) \in T_b B$ for
all $b\in B$. Hence $D\varphi(b)\in \b$, $D'\varphi \in \b$  for each $b\in B$.
This implies that ${\r}D\varphi = -D\varphi$, ${\r}D'\varphi = -D'\varphi$ and
$$
2\{ \varphi, \psi\}(b)= -\< D'\varphi(b), D' \psi(b)\> +
\< D\varphi(b), D \psi(b)\>
$$
But $D'\varphi(b)= Ad_{b^{-1}} D\varphi(b)$, $D'\psi(b)= Ad_{b^{-1}} D\psi(b)$
and $\< , \>$  is $Ad$-invariant. \qed

\medskip
\no For the next corollary note that $T^*_b(B)$ is a quotient of $T^*_b(G)$.

\begin{cor}
Let $\pi$ be the skew-symmetric 2-tensor on $T^*(G)$ corresponding to $w$. Pick
$b\in B$ and $\al,\be\in T^*_b(G)$. Then $\pi|b$ depends only  on the images of
$\al$ and $\be$ in $T^*_b(B)$.
\end{cor}

We will continue to use $\pi$ for the skew-symmetric 2-tensor on $T^*(B)$
induced by $\pi$ above.

\begin{rem}
An argument identical to that above proves that $w(k)$ is tangent to $K$. Hence
$\{ \cdot , \cdot \}$ induces a Poisson structure on $K$. With the above
structures $K$ and $B$ are sub Poisson Lie subgroups of the Poisson Lie group
$G$.
\end{rem}

We will need a formula for the Poisson tensor $\pi$ on $B$. We will use $\< ,
\>$ to identify $T^*(G)$ and $T(G)$. Under this identification $T^*_b(B)$ is
identified to $T_b(G)/T_b(B)$. We will identify this quotient with $dR_b \k$. We
let $\check{\pi}|_b$ denote the resulting
skew-symmetric 2-tensor on  $dR_b\k$. Finally we define $\check{\pi}^r\in C^\8
(B, (\La^2 \k)^*)$ by
$$
\check{\pi}|_b^r (x,y)= \check{\pi}|_b(dR_b x, dR_b y), \quad x,y\in\k
$$
We now recover formulae (2.25) of \cite{FR} (or \cite[Definition 4.2]{LR}) for
$\check{\pi}^r$.

\begin{lem}
\label{2.4}
$\check{\pi}|_b^r (x,y)= \< \rho_\k (Ad_{b^{-1}} x), \rho_\b (Ad_{b^{-1}} y)\>$.
\end{lem}
\proof Choose $\varphi, \psi\in C^\8 (G)$ with $\nab \varphi (b)= dR_b x$ and
$\nab \psi (b)= dR_b y$. Then
\begin{eqnarray*}
\check{\pi}|_b^r (x,y) & = & \check{\pi}|_b(dR_b x, dR_b y) \: =\:
\check{\pi}|_b (\nab \varphi (b),
 \nab \psi (b)) \\
& = & \pi|_b (d\varphi (b), d\psi(b)) \: =\: \{ \varphi, \psi\}(b) \\
& = & \half \< {\r} dL_b^{-1} \nab \varphi (b), dL_b^{-1} \nab \psi(b)\> -
\half \< {\r} dR_b^{-1} \nab \varphi (b), dR_b^{-1} \nab \psi(b)\> \\
& = & \half \< {\r} Ad_{b^{-1}} x, Ad_{b^{-1}} y\> - \half \< {\r}x,
y\>.
\end{eqnarray*}
But $x,y\in \k$ implies $\< {\r}x,y\> =0$. Hence,
\begin{eqnarray*}
\check{\pi}|_b^r (x,y) & = & \half \< \rho_\k (Ad_{b^{-1}} x), Ad_{b^{-1}} y\>-
  \half \< \rho_\b (Ad_{b^{-1}} x), Ad_{b^{-1}} y\> \\
  & = & \half \< \rho_\k (Ad_{b^{-1}} x), \rho_\b Ad_{b^{-1}} y\>-
  \half \< \rho_\b (Ad_{b^{-1}} x), \rho_\k Ad_{b^{-1}} y\> \\
  & = &\< \rho_\k (Ad_{b^{-1}} x), \rho_\b (Ad_{b^{-1}} y)\>.
\end{eqnarray*}
The last equality holds by skew-symmetry, see \cite[Lemma 4.3]{LR}. \qed

\medskip
\no We will abuse notation and drop the $\check{~}$ and $r$ in the notation for
$\check{\pi}^r$ henceforth.
\begin{rem}
The Poisson tensor on $K$, $\pi_K$, induced from the Skylanin bracket on $G$ is
the negative of the usual Poisson tensor on $K$ (see \cite{FR}, \cite{Lu}).  Throughout this
paper we let $\pi_K(k) = dL_k \, X\wedge Y - dR_k \, X\wedge Y$, where $X=\half \left(
\begin{smallmatrix}
0 & 1 \\ -1 & 0 \end{smallmatrix} \right) $ and $Y= \half \left(
\begin{smallmatrix} 0 & i \\ i & 0 \end{smallmatrix} \right)$.
\end{rem}

We now give  $G^n$ the product Poisson structure, hence $B^n$ inherits the
product structure.
We introduce more notation to deal with the product.
We let $\g_i \subset \g^n= \g \oplus ... \oplus \g$ be the image of $\g$ under
the embedding
into $i$-th summand. For  $\varphi\in C^\8 (G^n)$ we define
$$
D_i \varphi : G^n \to \g_i , \quad D_i' \varphi : G^n \to \g_i
$$
as follows. Let $g=(g_1,...,g_n)\in G^n$ and $\nu \in \g_i$, then
$$
\<D_i \varphi,\nu\>  =
\frac{d}{dt}|_{t=0} \varphi(g_1,..., e^{t\nu}g_i, ..., g_n)
$$
$$
\<D_i '\varphi,\nu\>  =
\frac{d}{dt}|_{t=0} \varphi(g_1,..., g_i e^{t\nu}, ..., g_n)
$$
Here we extend $\< , \>$ to $\g^n$ by
$$
\< \del, \ga\> =\sum_{i=1}^n \< \del_i, \ga_i\>
$$
for $\del=(\del_1,...,\del_n), \ga=(\ga_1,...,\ga_n)$. We define $d_i,\nab_i$ in
an analogous fashion. Finally define the Poisson bracket on $C^\8 (G^n)$ by
$$
\{ \varphi, \psi\}(g)= \half \sum_{i=1}^n [\< {\r} D_i' \varphi(g),
D_i'\psi(g)\> - \< {\r} D_i \varphi(g), D_i\psi(g)\>].
$$
As expected we obtain an induced Poisson bracket on $C^\8 (B^n)$ using the above
formula with $g\in G$ replaced by $b\in B$.

\medskip
Now let $\pi$ be the Poisson tensor on $G^n$ corresponding to the above Poisson bracket.  Let $\pi^{\#} \in Hom(T^*(G^n), T(G^n))$ be defined by $\be(\pi^\# (\al)) = \pi(\al , \be)$. Let $\varphi\in C^\8 (G^n)$. We have

\begin{defn}
The Hamiltonian vector field associated to $\varphi$ is the vector
field $X_\varphi 
\in C^\8 (G^n,  T(G^n))$ given by
$$
X_\varphi =\pi^\# d\varphi
$$
\end{defn}

We will need a formula for $X_\varphi$.

\begin{lem}
\label{ham}
Let $g=(g_1,...,g_n)$. Then $X_{\varphi}(g)= (X_1(g),..., X_n(g))$
where
$$
X_i(g)= \half[dR_{g_i} {\r} D_i \varphi(g) - dL_{g_i} {\r} D_i' \varphi(g)].
$$
\end{lem}
\proof We will use the formula
$$
\{\varphi, \psi\} = -d\psi(X_\varphi) = -\< X_\varphi , \nab \psi\>
$$
Here $\nab$ is the gradient with respect to $\< ,\>$ on $\g^n$, see above, 
hence $\nab\psi\!=\! ( \nab_1\psi,..., \nab_n\psi)$. We have

\begin{eqnarray*}
 \{\varphi, \psi\}(g) & = & \half \sum_{i=1}^n [\< \r D'_i \varphi(g), D'_i
\psi(g)\> -
 \< \r D_i \varphi(g), D_i \psi(g)\>] \\
 & = & \half \sum_{i=1}^n [\< \r D'_i \varphi(g), d L_{g_i}^{-1} \nab_i
\psi(g)\> -
 \< \r D_i \varphi(g), d R_{g_i}^{-1} \nab_i \psi(g)\>] \\
 & = & \half \sum_{i=1}^n [\< d L_{g_i} \r D'_i \varphi(g),  \nab_i \psi(g)\> -
 \< d R_{g_i} \r D_i \varphi(g),  \nab_i \psi(g)\>] \\
 & = & \half \< (X_1(g),..., X_n(g)), (\nab_1 \psi(g),..., \nab_n\psi(g))\> \\
 & = & \half \< X_\varphi (g), \nab\psi (g)\>
\end{eqnarray*}
\qed

\begin{rem}
Since $w(b)$ is tangent to $B^n$ the field $X_\varphi (b)$ will also be tangent
to $B^n$.
\end{rem}

\subsection{The dressing action of $K$ on $B^n$ and the action on $n$-gons in
$G/K$}
\label{dress} \label{S4.2}

The basic reference for this section is \cite{FR}. In that paper the authors
take $n=2$ and write $G=KB$. We will leave to the reader the task of comparing
our formulae with theirs.

In what follows we let  $G=SL_2(\C), \, K=SU(2)$, and $B$ be the
subgroup of $G$ consisting of upper-triangular matrices with
positive diagonal entries.  We let $\rho_B, \rho_K$ be the
projections relative to the decomposition $G= BK$. For the next
theorem (in the case $n=2$) see \cite[Formula 2.15]{FR}.

\begin{thm}
There is a Poisson action of $K$ on  the Poisson manifold $B^n$ given by
$$
k\cdot (b_1,...,b_n)= (b_1',...,b_n')
$$
with $b_i'= \rho_B(\rho_K(kb_1 \cdots b_{i-1}) b_i), \; 1\leq i\leq n$.
\end{thm}

\begin{defn}
The above action is called the {\bf dressing action} of $K$ on $B^n$.
\end{defn}
\begin{defn}
For $n=1$, we denote by $B_r=Kb$ the dressing orbit of $b$, where $b \in B$ and
$d(b \,*,*)=r$.
\end{defn}
We will also need the formula for the infinitesimal dressing action of $\k$ on
$B^n$. This action is given for $x\in \k$ by
$$
x\cdot (b_1,\cdots , b_n)= (\xi_1,...,\xi_n)\in T_b (B^n)
$$
with $\xi_i = dL_{b_i} \rho_{\b} Ad_{b^{-1}_i} \rho_\k Ad_{(b_1\cdots
b_{i-1})^{-1}} x$.
Note that $\xi_i \in T_{b_i} (B)$.

\begin{rem}
In order to pass from the $K$-action to the $\k$-action observe that
$\rho_K(bg)=\rho_K(g)$ and
$\rho_B(bg)=b\rho_B(g)$. Accordingly we may rewrite
the  $K$-action as  $k\cdot (b_1,...,b_n)= (b_1',...,b_n')$ with
$$
b_i'= b_i \rho_B(b^{-1}_i
\rho_K((b_1 \cdots b_{i-1})^{-1} k b_1 \cdots b_{i-1} ) b_i ), \; 1\leq i\leq n.
$$
\end{rem}

Recall, $*\in\h^3$ is the element fixed by the action of $K$,
$K\cdot*=*$.  Since $B$ acts simply-transitively on $G/K$ we have 

\begin{lem}
(i) The map $\Phi: B^n \to Pol_n(*)$ given by
$$
\Phi(b_1,...,b_n)= (*, b_1*, ..., b_1 \cdots b_n *)
$$
is a diffeomorphism.

(ii) The map $\Phi$ induces a diffeomorphism from $\{b\in B^n : b_1 \cdots b_n
=1\}$ onto $CPol_n(*)$.
\end{lem}

\no We now have

\begin{lem}
$\Phi$ is a $K$-equivariant diffeomorphism where $K$ acts on $B^n$ by the
dressing action
and on $Pol_n(*)$ by the natural (diagonal) action.
\end{lem}
\proof Let $k\cdot (b_1,...,b_n)= (b_1'',...,b_n'')$ be the pull-back to $B^n$
of the action
of $K$ on $Pol_n(*)$. Then
$$
b_1'' * = kb_1 *
$$
$$
b_1'' b_2'' * = kb_1 b_2*
$$
\centerline{$\vdots$}
$$
b_1'' \cdots b_n'' * = kb_1 \cdots b_n *
$$
We obtain
$$
b_1'' \cdots b_i''= \rho_B(k b_1 \cdots b_i)
$$
$$
b_i''= \rho_B((b_1''\cdots b''_{i-1})^{-1} k ( b_1 \cdots b_i))=
\rho_B((\rho_B (b_1''\cdots b''_{i-1}))^{-1} k  b_1 \cdots b_{i-1} b_i)
$$
$$
= \rho_B((\rho_B (kb_1\cdots b_{i-1}))^{-1} k  b_1 \cdots b_{i-1} b_i)=
\rho_B(\rho_K(k b_1 \cdots b_{i-1}) b_i)
$$
\qed

There is another formula for the dressing action of $K$ on $B^n$ that will be
useful.

\begin{lem}
With the above notation
$$
b_i'= \rho_B (\underbrace{\rho_K(\cdots (\rho_K}_{i-1~\hbox{times}} (k b_1) b_2
)\cdots b_{i-1}) b_i)
$$
\end{lem}
\proof Induction on $i$. \qed

\medskip
\no We obtain a corresponding formula  for the infinitesimal dressing action of
$\k$ on  $B^n$.

\begin{lem}
$x\cdot (b_1,...,b_n)= (\xi_1,...,\xi_n)$ where
$$
\xi_i = dL_{b_i} \rho_\b Ad_{b_i^{-1}} \rho_\k Ad _{b_{i-1}^{-1}}\cdots
\rho_\k Ad_{b_1^{-1}}x.
$$
\end{lem}

\no We now draw an important consequence.

\begin{lem}
The map $\Phi$ induces a diffeomorphism between $B_r^n = B_{r_1}
\times \cdots \times B_{r_n}$ and the configuration space of open based $n$-gon
linkages
$\t{N}_r$, where if $b \in B_r^n$, then $r= (r_1,..,r_n)$ and $d(b_{1}\cdots b_i
*,b_1 \cdots
b_{i-1} *)=r_i$, for all $1\leq i \leq n$.
\end{lem}
\proof Let $b\in B_r^n$ be given.
Then $\Phi(b)= (*, b_1 *,..., b_1 \cdots b_n *)$.
The $K$-orbit of $\Phi(b)$ is
$$
[ *, kb_1 *,..., k(b_1 \cdots b_n)*]
$$
The $i$-th edge $e_i$ of $\Phi(b)$ is the geodesic
segment joining $k b_1 \cdots b_{i-1}*$ to
$k b_1 \cdots b_i *$. Clearly this is congruent (by $k b_1 \cdots b_{i-1}$) to
the segment connecting $*$ to $b_i *$. \qed

\begin{cor}
The symplectic leaves of $B^n$ map to the configuration spaces  $\t{N}_r$ under
$\Phi$.
\end{cor}
\proof The $B_r^n$ are the symplectic leaves of $B^n$. \qed

\subsection{The moduli space $M_r$ as a symplectic quotient} \label{S4.3}

We have seen that $\Phi$ induces a diffeomorphism from
$\{ (b_1,..., b_n)\in B_r^n : b_1 \cdots b_n =1\}/K$
to the moduli space $M_r$ of closed $n$-gon linkages
in $\h^3$ modulo isometry. In this section we will prove that the map $\varphi:
B^n\to B$ given by
$\varphi(b_1 ,\cdots, b_n)= b_1 \cdots b_n$ is a momentum map for the (dressing)
$K$ action on $B^n$. Hence, $M_r$ is a symplectic quotient, in particular
it is a symplectic manifold if $1$ is a regular value by Lemma \ref{4.3} below.

The definition of a momentum map for a Poisson action of a Poisson Lie group was
given in \cite{Lu}.

\begin{defn}
Suppose that $K$ is a Poisson Lie group, $(M,\pi)$ is a Poisson manifold, and
$K\times M \to M$ is a Poisson action. Let $x\in \k$, $\al_x$ be the extension
of $x\in\k= (\k^*)^*$ to a
right-invariant 1-form on $K^*$, and $\hat{x}$
be the induced vector field on $M$. Then a map $\varphi: M\to K^*$ is a {\bf
momentum map} if it satisfies the equation
$$
- \pi^\# \varphi^* \al_x = \hat{x}
$$
\end{defn}

\begin{rem}
For the definition of $K^*$, the {\bf dual Poisson Lie group}, see
\cite{Lu}. In our case $K^*=B$. 
\end{rem}

\no The next lemma is proved in \cite{Lu} and \cite{Lu3}.  We include a
proof here for completeness.

\begin{lem}
\label{4.299} Suppose that $M$ is a symplectic manifold, $K$ is a
Poisson Lie group and $K\times M\to M$ is a Poisson action with an
equivariant momentum map $\varphi: M\to K^*$. Assume $1$ is a
regular value of $\varphi$. Then $\varphi^{-1}(1)/K$ is a
symplectic orbifold with the symplectic structure given by taking
restriction and quotient of the symplectic structure on $M$.  If
we assume further that the isotropy subgroups of all $x \in
\varphi^{-1}(1)$ are trivial then $\varphi^{-1}(1)/K$ is a
manifold.
\end{lem}
\proof Let $\om$ be the symplectic form on $M$ and $m\in \varphi^{-1}(1)\subset
M$. Let $V_m$ be a
complement to $T_m \varphi^{-1}(1)$ in $T_m M$. Let $\k \cdot m\subset T_m M$
be the tangent space to the orbit $K\cdot m$. Hence $\k_m = \{ \hat{x}(m) : x\in
\k\}$.
 We first prove the identity
$$
\om_m \hat{x}(m) = -\varphi^* \al_x |_m
$$
Here we use $\om_m$ to denote the induced map $T_m M \to T^*_m M$ as well as the
symplectic
form evaluated at $m$. Indeed we have the identity
$$
-\pi^\# \varphi ^* \al_x = \hat x
$$
Applying $\om$ we get
$$
\om \hat x = -\varphi^* \al_x
$$
We claim that if $m\in \varphi^{-1}(1)$ the $\k \cdot m$ is orthogonal
 (under $\om_m$) to $T_m(\varphi^{-1}(1))$, in particular it is
totally-isotropic. Let $\hat{x}(m)\in \k\cdot m$ and
$u\in T_m (\varphi^{-1}(1))$. Then
$$
\om_m (\hat{x}(m),u)= -\varphi^*\al_x|_m (u)= -\al_x|_m (d\varphi_m u)
$$
But $d\varphi_m u=0$ and the claim is proved. Hence, the restriction of
$\om_m$ to $T_m(\varphi^{-1}(1))$ descends to $T_m(\varphi^{-1}(1))/\k\cdot m=
T_{K\cdot m} (\varphi^{-1}(1)/K)$.
We now prove that the induced form is nondegenerate.

To this end we claim that $\k \cdot m$ and $V_m$ are dually paired by $\om_m$.
We draw two conclusions from the
hypothesis that $1$ is a regular value for $\varphi$. First by \cite[Lemma
4.2]{FR} the map $\k\to \k \cdot m$
given by $x\mapsto \hat{x}(m)$ is an isomorphism. Second, $d\varphi_m:V_m \to
\k^*$ is an isomorphism.
Let $\{x_1,...,x_N\}$ be a basis
for $\k$, whence $\{\hat{x}_1(m),..., \hat{x}_N(m)\}$ is a basis for $\k \cdot
m$. We want to find a
basis $\{v_1,...,v_N\}$ for $V_m$ so that $\om_m(\hat{x}_i(m),v_j)= \del_{ij}$.
Choose a basis $\{v_1,...,v_N\}$ for $V_m$ such that
$\{d\varphi_m v_1,...,d\varphi_m v_N\}\subset \k^*$ is dual to
$\{x_1,...,x_N\}$. Then
$$
\om_m(\hat{x}_i(m),v_j)= -\varphi^* \al_{x_i} |_m (v_j)= -\al_{x_i}|_m (d\varphi
v_j)=
-d\varphi v_j(x_i)= -\del_{ij}
$$
As a consequence of the previous claim, the restriction of $\om_m$ to $\k\cdot
m\oplus V_m$ is
nondegenerate. We then have the orthogonal complement $(\k\cdot m\oplus
V_m)^{\perp}$ is a complement to
$\k\cdot m\oplus V_m$ and $\om_m| (\k\cdot m\oplus V_m)^{\perp}$ is
nondegenerate. But then
$(\k\cdot m\oplus V_m)^{\perp}$ maps isomorphically to $T_{K\cdot m}
(\varphi^{-1}(1)/K)$. \qed

\medskip
\no We will also need
\begin{lem}
\label{reg}
$\varphi(m)$ is a regular value for $\varphi$ if and only if $\k_m = \{x \in
\k: \hat{x}(m)=0\} = 0$.
\end{lem}
\proof  Let $x\in \k$. Then $x\in (Im\,d\varphi|_m)^\bot \Leftrightarrow
\varphi^* \al_x = 0 \Leftrightarrow 0 = -\pi \varphi^* \al_x = \hat{x}(m)$.
\qed

\medskip
We now begin the proof that $\varphi$ is a momentum map for the dressing action
of $K$ on $B^n$.
We will need some notation. Let $x\in \k = \b^*$. Recall that $\al_x$ is the
extension of $x$ to a
right-invariant 1-form on $B$.
Thus if $\zeta\in T_b (B)$ we have
$$
\al_x |_b (\zeta)= \< dR_b x, \zeta\>.
$$

\begin{lem}
\label{4.3}
$-\pi^\# \al_x |_b = dL_b \rho_\b Ad_{b^{-1}} x$.
\end{lem}
\proof Let $y\in \k$. It suffices to prove that
$$
\pi|_b (\al_x(b), \al_y(b))= -\al_y |_b (dL_b \rho_{\b} Ad_{b^{-1}}x )
$$
Now
$$
\pi|_b (\al_x(b), \al_y(b))= \pi|_b (dR_b x, dR_b y)= \pi^r|_b (x,y)= \<
\rho_{\k}(Ad_{b^{-1}} x), \rho_{\b}(Ad_{b^{-1}}y)\>
$$
according to Lemma \ref{2.4}. Also
$$
-\al_y |_b  (dL_b \rho_\b Ad_{b^{-1}} x)= -\< dR_b y, dL_b \rho_\b Ad_{b^{-1}}
x\>
$$
$$
= -\< Ad_{b^{-1}} y, \rho_\b Ad_{b^{-1}} x\>=  -\< \rho_\k Ad_{b^{-1}} y,
\rho_\b Ad_{b^{-1}} x\>= \< \rho_\k Ad_{b^{-1}} x, \rho_\b Ad_{b^{-1}} y\>
$$
\qed

\begin{lem}
\ \ $\varphi^* \al_x|_b = (\al_{x_1} |_{b_1}, \cdots, \al_{x_n} |_{b_n})$, 
\ where \ $x_1=x$ \ and
$x_i= \rho_\k (Ad_{(b_1\cdots b_{i-1})^{-1}} x)$, $2\leq i\leq n$.
\end{lem}
\proof We will use the following formula (the product rule).
Let $b=(b_1,...,b_n)$ and $\zeta= (\zeta_1,...,\zeta_n)\in T_b(B^n)$. Then
$$
d\varphi_b(\zeta)= (dR_{b_2\cdots b_n} \zeta_1+ dL_{b_1} dR_{b_3\cdots b_n}
\zeta_2+\cdots + dL_{b_1\cdots b_{n-1}}\zeta_n)
$$
Hence
\begin{eqnarray*}
  (\varphi^* \al_x)|_b(\zeta) & = & \al_x|_{b_1 \cdots b_n} (dR_{b_2\cdots
b_n}\zeta_1+ \ dL_{b_1} dR_{b_3\cdots b_n} \zeta_2+ \
  \cdots + \ dL_{b_1\cdots b_{n-1}}\zeta_n) \\
  & = & \< dR_{b_1\cdots b_n} x, dR_{b_2\cdots b_n} \zeta_1\> + \< dR_{b_1\cdots
b_n} x,
  dL_{b_1}dR_{b_3\cdots b_n} \zeta_2\>+\cdots \\
  &   & \: \: \: +\< dR_{b_1\cdots b_{n}}x , dL_{b_1\cdots b_{n-1}}\zeta_n\>
  \\
  & = & \< dR_{b_1} x, \zeta_1\> + \< dR_{b_2} Ad_{b_1^{-1}} x, \zeta_2\> +
\cdots +
  \< dR_{b_n} Ad_{(b_1\cdots b_{n-1})^{-1}} x, \zeta_n\> \\
  & = & \sum_{i=1}^n \al_{x_i}(\zeta_i)
\end{eqnarray*}

\begin{prop}
$\varphi$ is an equivariant momentum map for the dressing action of $K$ on
$B^n$.
\end{prop}
\proof To show that $\varphi$ is a momentum map we have to check that
$$
- \pi^\# \varphi^* \al_x |_b = x\cdot (b_1,...,b_n)
$$
But $- \pi^\# \varphi^* \al_x= (-\pi^\# \al_{x_1}, \cdots, -\pi^\# \al_{x_n})$ and the
result follows from the previous two lemmas. To show that $\varphi$ is
equivariant we have to check that $\varphi(k\cdot (b_1,...,b_n))= k\cdot b_1 ...
b_n$.
This is obvious from the point of view of polygons.
\qed

\medskip
\no As a consequence of the above proposition we obtain

\begin{thm}
\label{sq}
The map $\Phi$ carries the symplectic quotient $(\varphi |_{B_r^n})^{-1}(1)/K$
diffeomorphically to the moduli space of $n$-gon linkages $M_r$.
\end{thm}
\begin{rem} We obtain a symplectic structure on $M_r$ by transport of structure.
\end{rem}

\subsection{The bending Hamiltonians}

In this section we will compute the Hamiltonian vector fields $X_{f_j}$ of
the
functions
$$
f_j(b)= tr((b_1 \cdots b_j)(b_1 \cdots b_j)^*), \; 1\leq j\leq n.
$$
Throughout the rest of the paper, we will assume $G=SL_2(\C)$. Then $G=BK$,
where
$B =\{\left(
\begin{smallmatrix}
a & z \\
0 & \; a^{-1}
\end{smallmatrix}
\right) \in SL_2(\C) | a\in \R_+, \; z \in  \C \}$ and $K=SU(2)$.

We will use the following notation. If $A\in M_m(\C)$ then
$A^0= A- \frac{1}{m} tr(A)I$ will be its projection to the traceless matrices.

\begin{thm}
\label{tree} Define $F_j: B^n \to \k$ for $b = (b_1, b_2, ..., b_n)$ by 
$$ F_j(b)= \sqrt{-1}[(b_1 \cdots b_j)(b_1 \cdots b_j)^*]^0 
$$
Then $X_{f_j}(b) =
(F_j(b)\cdot (b_1,...,b_j),0, ...,0)$ where $\cdot$ is the
infinitesimal dressing action of $\k$ on $B^j$, see \S
\ref{dress}.
\end{thm}
\proof It will be convenient to work on $G^n$ and then restrict to $B^n$.
By the formula for $X_\varphi$ of Lemma \ref{ham} it suffices to compute $D_i
f_j$
and $D_i' f_j$. We recall that
$$
D'_i \varphi (g)= Ad_{g_i^{-1}} D_i \varphi(g)
$$
hence it suffices to compute $D_i f_j(g)$. We first reduce to computing $D_1
f_j$ by

\begin{lem}
\label{addi}
$D_i f_j(g)= Ad_{(g_1\cdots g_{i-1})^{-1}} D_1 f_j(g)$.
\end{lem}
\proof By definition
$$
\< D_i f_j(g),\nu\>= \frac{d}{dt}|_{t=0} f_j(g_1,..., e^{t\nu}g_i,..., g_n)
$$
But it is elementary that
$$
f_j(g_1,..., e^{t\nu}g_i,..., g_n)= f_j((Ad_{g_1\cdots g_{i-1}}
e^{t\nu})g_1,...,g_n)
$$
Differentiating at $t=0$ we obtain
$$
\< D_i f_j(g),\nu\>= \< D_1 f_j(g), Ad_{g_1\cdots g_{i-1}}\nu\> =
\< Ad_{(g_1\cdots g_{i-1})^{-1}} D_1 f_j(g),\nu\>
$$
\qed

\no We next have

\begin{lem}
$D_1 f_j(g)= F_j(g)$.
\end{lem}
\proof By definition
\begin{eqnarray*}
  \< D_1 f_j(g),\nu\> & = & \frac{d}{dt}|_{t=0} tr[(e^{t\nu}g_1\cdots
g_j)(e^{t\nu}g_1\cdots g_j)^*] \\
  & = & tr[(\nu g_1\cdots g_j)(g_1\cdots g_j)^* + (g_1\cdots g_j)(\nu g_1\cdots
g_j)^*] \\
  & = & tr[(\nu g_1\cdots g_j)(g_1\cdots g_j)^*] + tr[(\nu g_1\cdots
g_j)(g_1\cdots g_j)^*]^* \\
  & = & tr[(\nu g_1\cdots g_j)(g_1\cdots g_j)^*] + \overline{tr[(\nu g_1\cdots
g_j)(g_1\cdots g_j)^*]} \\
  & = & 2 Re tr [(\nu g_1\cdots g_j)(g_1\cdots g_j)^*] \\
  & = & 2Im \sqrt{-1} tr [\nu (g_1\cdots g_j)(g_1\cdots g_j)^*]
\end{eqnarray*}

\no Since $\nu\in sl_2(\C)$ we may replace
$$(g_1\cdots g_j)(g_1\cdots g_j)^*$$
by its trace\-less projection
$$[(g_1\cdots g_j)(g_1\cdots g_j)^*]^0 \;.$$

\no Since $tr$ is complex bilinear we obtain
\begin{eqnarray*}
\< D_1 f_j(g),\nu\> & = & 2Im tr (\nu \sqrt{-1} [ (g_1\cdots g_j)(g_1\cdots g_j)
  ^*]^0) \\
  & = & \< \sqrt{-1} [ (g_1\cdots g_j)(g_1\cdots g_j)^*]^0,\nu\>.
\end{eqnarray*}
\qed

Now we restrict to $B^n$ and substitute into our formula for $X_{f_j}(b)$ in
Lemma \ref{ham}.
We obtain

\begin{lem}
$ (X_{f_j})_i= (D_1 f_j(b)\cdot (b_1,..., b_n))_i$= the $i$-th component of the
infinitesimal dressing action of $D_1 f_j(b)\in \k$.
\end{lem}
\proof By Lemma \ref{ham} we have
$$
(X_{f_j})_i= \half [dR_{b_i} \r D_i f_j(b)- dL_{b_i} \r D_i ' f_j(b)]
$$
$$
= \half [dR_{b_i} \r Ad_{(b_1\cdots b_{i-1})^{-1}} D_1 f_j(b)-
dL_{b_i} \r Ad_{b_i^{-1}} Ad_{(b_1\cdots b_{i-1})^{-1}} D_1 f_j(b)]
$$
We  write
$$
Ad_{(b_1\cdots b_{i-1})^{-1}} D_1 f_j(b)= X_1 +\eta_1
$$
with $X_1\in \k$ and $\eta\in \b$. Hence
$$
X_1 = \rho_\k (Ad_{(b_1\cdots b_{i-1})^{-1}}D_1 f_j(b))
$$
$$
\eta_1= \rho_\b (Ad_{(b_1\cdots b_{i-1})^{-1}}D_1 f_j(b))
$$
Then
$$
\r Ad_{(b_1\cdots b_{i-1})^{-1}}D_1 f_j(b) =X_1 - \eta_1
$$
We write
$$
Ad_{b_i^{-1}}X_1 = X_2 +\eta_2, X_2\in \k, \eta_2\in \b
$$
Then $X_2= \rho_\k ( Ad_{b_i^{-1}}X_1) , \eta_2= \rho_\b ( Ad_{b_i^{-1}}X_1)$
and
$$
Ad_{b_i^{-1}}(X_1 +\eta_1) = X_2 +\eta_2 + Ad_{b_i^{-1}}\eta_1
$$
Hence
\begin{eqnarray*}
  \r Ad_{b_i^{-1}} Ad_{(b_1\cdots b_{i-1})^{-1}}D_1 f_j(b) & = & X_2 - \eta_2 -
Ad_{b_i^{-1}} \eta_1
  \\
  & = & X_2 +\eta_2  -2\eta_2 - Ad_{b_i^{-1}}\eta_1 \\
  & = & Ad_{b_i^{-1}} X_1 - 2\eta_2 - Ad_{b_i^{-1}}\eta_1
\end{eqnarray*}
Hence
$$
(X_{f_j})_i= \half [ dR_{b_i} X_1 - dR_{b_i}\eta_1 - dL_{b_i} Ad_{b_i^{-1}} X_1
+ 2 dL_{b_i} \eta_2 + dL_{b_i} Ad_{b_i^{-1}} \eta_1]
$$
But $dL_{b_i} Ad_{b_i^{-1}}= dR_{b_i}$ and we obtain
$$
(X_{f_j})_i= dL_{b_i}\eta_2
$$
Since $\eta_2= \rho_\b ( Ad_{b_i^{-1}} \rho_\k ( Ad_{(b_1\cdots
b_{i-1})^{-1}}D_1 f_j(b)))$ the
lemma follows.
\qed

\medskip
\no With this Theorem \ref{tree} is proved.
\qed

\subsection{Commuting Hamiltonians}

In this section we will show the functions $$ f_j(b)= tr((b_1
\cdots b_j)(b_1 \cdots b_j)^*), \; 1 \leq j \leq n $$ Poisson
commute.  The proof is due to Hermann Flaschka.

\begin{prop}
$\{f_j,f_k \} = 0$  for all $j,k$.
\end{prop}

\proof  Again we will work on $G^n$ and then restrict to $B^n$.  Without loss of
generality we let $j \leq k$.

\no Recall from Lemma \ref{addi}
$$
D_i f_j(g)= Ad_{(g_1\cdots g_{i-1})^{-1}} D_1 f_j(g).
$$
It is easily seen that
$$
D_i f_j(g) = D_{i-1}' f_j(g), \;\textnormal{for}\; 1 \leq i \leq j.
$$

\no We now have,
\begin{eqnarray*}
\{ f_j, f_k\}(g) & = & \half \sum_{i=1}^n [\< {\r} D_i' f_j(g), D_i'f_k(g)\> -
\< {\r} D_i f_j(g), D_i
f_k(g)\>] \\
  & = & \half \sum_{i=1}^j [\< {\r} D_i' f_j(g), D_i' f_k(g)\> - \< {\r} D_i
f_j(g), D_i f_k(g)\>] \\
  & = & \half [\< {\r} D_j' f_j(g), D_i' f_k(g)\> - \< {\r} D_1 f_j(g), D_1
f_k(g)\>] \\
  & = & \half [\< {\r} Ad_{(g_1\cdots g_j)^{-1}} D_1 f_j(g), Ad_{(g_1\cdots
g_j)^{-1}} D_1 f_k(g)\> - \< {\r} D_1 f_j(g), D_1 f_k(g)\>]
  \\
  & = & \half \< {\r} Ad_{(g_1\cdots g_j)^{-1}} D_1 f_j(g), Ad_{(g_1\cdots
g_j)^{-1}} D_1 f_k(g)\>
\end{eqnarray*}

\no since $D_1 f_i(g) \in \k$ for all $i$.
\
The proposition follows if we can show
$$
\< {\r} Ad_{(g_1\cdots g_j)^{-1}} D_1 f_j(g), Ad_{(g_1\cdots g_j)^{-1}} D_1
f_k(g)\> = 0.
$$

\no It follows from the proof of Theorem \ref{tree} that

\begin{eqnarray*}
Ad_{(g_1\cdots g_j)^{-1}} D_1 f_j(g) & = & \sqrt{-1} Ad_{(g_1\cdots g_j)^{-1}}[
(g_1 \cdots g_j)(g_1 \cdots
g_j)^*]^0 \\
  & = & \sqrt{-1} [(g_1 \cdots g_j)^*(g_1 \cdots g_j)]^0 \in \k.
\end{eqnarray*}
Hence,
\begin{eqnarray*}
  \{ f_j, f_k \}(g) & = & \< {\r} Ad_{(g_1\cdots g_j)^{-1}} D_1 f_j(g),
Ad_{(g_1\cdots g_j)^{-1}} D_1 f_k(g)\> \\
  & = & \< Ad_{(g_1\cdots g_j)^{-1}} D_1 f_j(g), Ad_{(g_1\cdots g_j)^{-1}} D_1
f_k(g)\> \\
  & = & \< D_1 f_j(g), D_1 f_k(g)\> \\
  & = & 0
\end{eqnarray*}
since $\<,\>$ is $Ad$-invariant.  This proves the proposition on
$G^n$. The result then holds when we restrict to $B^n$. \qed

\subsection{The Hamiltonian flow}

In this section we compute the Hamiltonian flow, $\varphi_k^t$, associated
to
$f_k$.
\

Recall from Theorem \ref{tree} the Hamiltonian field for $f_k$ is
given by $X_{f_j}(b)= (F_j(b)\cdot (b_1,...,b_j), 0,...,0)$ where
$\cdot$ is the infinitesimal dressing action of $K$ on $B^n$.  We
now need to solve the system of ordinary differential equations

\begin{equation*}
(*)
 \begin{cases}
  \frac{db_i}{dt} = (F_j(b) \cdot (b_1,...,b_j))_i, \; 1 \leq i \leq j\\
  \frac{db_i}{dt} = 0, \; j+1 \leq i \leq n
 \end{cases}
\end{equation*}

\begin{lem}
$D_1 f_j(b)= F_j(b)$ is invariant along solution curves of (*).
\end{lem}

\proof It suffices to show $\varphi_j(b)=b_1\cdots b_j$ is constant along
solution curves.

\no Let $b(t) = (b_1 (t),..., b_n (t))$ be a solution of $X_{f_j}$.  Then

\begin{eqnarray*}
\frac{d}{dt} \varphi_j (b(t)) & = & \frac{db_1}{dt}(t) b_2(t) \cdots b_j(t) +
b_1(t) \frac{db_2}{dt}(t) \cdots b_j(t)
   + \cdots + b_1(t) b_2(t) \cdots \frac{db_j}{dt}(t) \\
   & = & \half [(\r D_1 f_j(b(t))) b_1(t) - b_1(t) \r D_1 ' f_j(b(t))] b_2(t)
\cdots b_j(t) \\
   &   &\; + b_1(t)  \half [(\r D_2 f_j(b(t))) b_2(t) - b_2(t) \r D_2 '
f_j(b(t))] b_3(t) \cdots b_j(t) + \cdots +\\
   &   &\; + b_1(t) b_2(t) \cdots  b_{j-1}(t) \half [(\r D_j f_j(b(t))) b_j(t) -
b_j(t) \r D_j ' f_j(b(t))] \\
   & = &\half [\r (D_1 f_j(b(t))) b_1(t) \cdots b_j(t) - b_1(t)\cdots b_j(t) \r
(D_j ' f_j(b(t)))] \\
   & = &\half [(D_1 f_j(b(t))) b_1(t) \cdots b_j(t) - b_1(t) \cdots b_j(t) (D_j
' f_j(b(t)))] \\
   & = &\half [(D_1 f_j(b(t))) b_1(t) \cdots b_j(t) - b_1(t) \cdots b_j(t)
(Ad_{(b_1\cdots b_j)^{-1}}D_1 f_j(b(t)))] \\
   & = &\half [(D_1 f_j(b(t))) b_1(t) \cdots b_j(t) - (D_1 f_j(b(t))) b_1(t)
\cdots b_j(t)] = 0
\end{eqnarray*}

\no Thus $\varphi (b)$ is constant along solution curves of $X_{f_j}$, proving
the lemma.
\qed

\begin{rem}
It also follows from the previous proof that $f_j(b)$ is constant along solution
curves of (*).
\end{rem}

Let $b = \left(
\begin{smallmatrix}
a &z \\
0 & a^{-1}
\end{smallmatrix}
\right)$
with $a \in \R_+$ and $z \in \C$,
then it follows from a simple calculation that
$$
det(F_1(b)) = \frac {1} {4}(a^4 +a^{-4}+ |z|^4 - 2 + 2a^2 |z|^2 +2a^{-2} |z|^2).
$$
Since $a > 0$ we see that $a^4 + a^{-4} \geq 2$ with equality if $a=1$.
Therefore, $det(F_1(b)) \geq 0$ with equality iff $b=1$.  From the above
argument it follows that $det(F_j(b)) \geq 0$ with equality iff $b_1
\cdots b_j =1$. 

It is also an easy calculation to show $$ det(F_j(b)) =
\textstyle{\frac{1}{4}} f_j(b)^2 -1,\; \ \forall b \in B^n $$
\begin{lem}
The curve $exp(tF_j(b))$is periodic with period $2\pi / \sqrt{\frac{1}{4} f_j(b)^2 -
1}$
\end{lem}

\proof To simplify notation, let $X=F_j(b) \in \k$. Then
$$
X^{-1} = -\frac {1} {det(X)} X
$$
giving us
$$
X^2 = -(det(X))X^{-1} X = -det(X) I
$$
So,
\begin{eqnarray*}
  \exp tX & = & \sum_{n=0}^\8 \frac {t^n X^n} {n!} \\
  & = & \sum_{n=1}^\8 \frac{(-1)^n (tdet(X))^n} {(2n)!}I + \sum_{n=1}^\8
\frac{(-1)^n (tdet(X))^n} {(2n+1)!} \frac
  {X}{\sqrt{det(X)}} \\
  & = & \cos \left( t \sqrt{det(X)} \right)I + \frac{\sin \left(t \sqrt{det(X)}\right)}{\sqrt{det(X)}}X
\\
  & = & \cos \left(t \sqrt{\textstyle{\frac{1}{4}} f_j(b)^2 -1}\right)I + \frac{\sin \left(t \sqrt{\frac{1}{4}
f_j(b)^2 -1}\right)}{\sqrt{\frac{1}{4} f_j(b)^2 -1}}F_j(b)
\end{eqnarray*}

\no Therefore the curve is periodic with period $2\pi/\sqrt{\frac{1}{4}
f_j(b)^2 -1}$. \qed

\medskip
\no We can now find a solution to the system (*)
\begin{prop}
\label{flow} Suppose $P \in M_r$ has vertices given by
$b_1,...,b_n$.  Then $P(t)=\varphi_k^t(P)$ has vertices given by
$b_1(t),...b_n(t)$ where $$ b_i(t) = (\exp (tF_j(b)) \cdot
(b_1,..., b_j))_i, \; 1 \leq i \leq k $$ $$ b_i(t)=b_i, \; k+1
\leq i \leq n. $$ Here $\cdot$ is the dressing action of $K$ on
$B^j$.
\end{prop}
\proof This follows from $F_j(b)$ being constant on solution curves of (*).  We
can see immediately that the $b_i$'s
are solutions curves of our system of ordinary differential equations.
\qed

\begin{cor}
The flow $\varphi_k^t(P)$ is periodic with period $2\pi / \sqrt{\frac{1}{4} f_j(b)^2
-1}$.
\end{cor}

\begin{rem}
If the k-th diagonal is degenerate ($b_1 \cdots b_k = 1$) then P
is a fixed point of $\varphi_k^t$.  In this case the flow has
infinite period.
\end{rem}

\no Let $\ell_k(b) = 2 \cosh^{-1} (\half f_k(b))$, then
$$
d\ell_k = \frac{1}{\sqrt{\frac{1}{4} f_k(b)^2 - 1}} df_k
$$
and consequently
$$
X_{\ell_k} = X_{f_k}/\sqrt{\textstyle{\frac{1}{4}} f_k^2 - 1}
$$
where $X_{\ell_k}$ is the Hamiltonian vector field associated to
$\ell_k$.  Since
$f_j$ is a constant of motion, $X_{\ell_k}$ is constant along solutions of (*)
as well.  Let
$\Psi_k^t$ be the flow of $X_{\ell_k}$.  We have the following

\begin{prop}
Suppose $P \in M'_r$ has vertices $b_1,...,b_n$.  Then $P(t) = \Psi_k^t(P)$ has
vertices $b_1(t),...,b_n(t)$ given by
$$
b_i(t) = \exp \Bigl( (tF_j(b))/\sqrt{\textstyle{\frac{1}{4}} f_k(b)^2 - 1}\Bigr)
\cdot (b_1,...,b_j))_i,
\; 1 \leq i \leq k
$$
$$
b_i(t)=b_i, \; k+1 \leq i \leq n
$$
where $\cdot$ is the dressing action of $K$ on $B^n$.
\end{prop}

Thus $\Psi_k^t$ is periodic with period $2\pi$ and rotates a part
of $P$ around the $k$-th diagonal with constant angular velocity 1
and leaves the other part fixed.

\subsection{Angle variables, the momentum polyhedron and a new proof of
involutivity}

We continue to assume that our n-gons are triangulated by the
diagonals $\{d_{1i},\; 3 \leq i \leq n-1 \}$.  We assume $P \in
M^o_r$ so none of the $n-2$ triangles, $\triangle_1,\triangle_2,
..., \triangle_{n-2}$, created by the above diagonals are
degenerate.  We construct a polyhedral surface $S$ bounded by $P$
by filling in the triangles $\triangle_1, \triangle_2, ...,
\triangle_{n-2}$.  Hence, $\triangle_1$ has edges $e_1, \,e_2,
\,$and$\,d_{13}$, \;$\triangle_2$ has edges $d_{13}, \,e_3,
\,$and$\,d_{14}$, ..., and $\triangle_{n-2}$ has edges $d_{1,n-1},
\,e_{n-1}, \,$and$\,e_n$.

We define $\hat{\th}_i$ to be the oriented dihedral angle measured
from $\triangle_i$ to $\triangle_{i+1}, \; 1 \leq i \leq n-3$.  We
define the $i$-th angle variable $\th_i$ by $$ \th_i = \pi -
\hat{\th}_i, \; 1 \leq i \leq n-3. $$

\begin{thm}
$\{\th_1,...,\th_{n-3}\}$ are angle variables, that is we have
\begin{itemize}
 \item [(i)]  $\{\ell_i, \th_j \} = \del_{ij}$
 \item [(ii)]  $\{\th_i, \th_j \} = 0$.
\end{itemize}
\end{thm}
\proof The proof is identical to that of \cite[\S4]{KM2}.
\qed

\medskip
We next describe the momentum polyhedron $B_r$ for the action of
the above (n-3)-torus by bendings.  Hence, $$ B_r=
\{\ell(M_r)\subset (\R_{\geq 0})^{n-3}: \; \ell
=(\ell_1,...,\ell_{n-3})\}. $$

Let $(\ell_1,...,\ell_{n-3}) \in (\R_{\geq 0})^{n-3}$ be given.
We first consider the problem of constructing the triangles,
$\triangle_1, \triangle_2, ..., \triangle_{n-2}$ above.  We note
that there are three triangle inequalities $E_i(\ell,r), \; 1 \leq
i \leq n-2$, among the $r_i$'s and $\ell_j$'s that give necessary
and sufficient conditions for the existence of $\triangle_i$. Once
we have obtained the triangles $\triangle_1, \triangle_2, ...,
\triangle_{n-2}$, we can glue them along the diagonals $d_{1i}, \;
3\leq i \leq n-1$, and obtain a polyhedron surface $S$ and a n-gon
$P$.  We obtain

\begin{thm}
The momentum polyhedron $B_r \subset (\R_{\geq0})^{n-3}$ is
defined by the $3(n-2)$ triangle inequalities
\begin{eqnarray*}
|r_1-r_2| \leq & \ell_1 & \leq r_1 +r_2 \\
|\ell_1-r_3| \leq & \ell_2 & \leq \ell_1 +r_3 \\
& \vdots & \\
|\ell_{n-4}-r_{n-2}| \leq & \ell_{n-3} & \leq \ell_{n-4}+r_{n-2} \\
|r_{n-1}-r_n| \leq & \ell_{n-3} & \leq r_{n-1}+r_n
\end{eqnarray*}
Here $r=(r_1,...,r_n)$ is fixed, the $\ell_i$'s, $1 \leq i \leq n-3$, are the
variables.
\end{thm}

\no As a consequence we have

\begin{thm}
The functions $\ell_1, \ell_2, ... , \ell_{n-3}$ on $M_r$ are functionally
independent.
\end{thm}

The theorem follows from Corollary 4.45.  We will apply the next
lemma with $M= \ell^{-1}(B_r^o)$, the inverse image of the interior of
the momentum polyhedron under $\ell = (\ell_1, \ell_2, ... ,
\ell_{n-3})$.  Then $M \simeq B_r^o \times (S^1)^{n-3}$. 

\begin{lem} \label{funind}
\ \ Suppose \ $M = M_r^o$ is a \ connected \ real-analytic \ manifold \ and 
$F=(f_1, ... ,f_k):M^n \to \R^k, \, n \geq k,$ is a real-analytic map
such that $F(M)$ contains a k-ball.  Then the 1-forms $df_1, ..., df_k$ are
linearly independent over $C^{\8}(M)$.
\end{lem}

\proof  Since the 1-forms $df_1, df_2, ..., df_k$ are real-analytic,
the set of points $x \in M$ such that $df_1|_x, ..., df_k|_x$ are not
independent over $\R$ is an analytic subset $W$ of $M$.  Let $M^0 = M
- W$. Hence either $M^o$ is empty or it is open and dense.  But by
Sard's Theorem, $F(W)$ has measure zero.  Since $F(M)$ does not have
measure zero, $M \ne W$ and $M^o$ is nonempty,  hence open and dense.
Therefore, if there exists $\varphi_1, ... , \varphi_k
\in C^{\8}(M)$ such that $\sum_{i=1}^k \varphi_i df_i = 0$ then
$\varphi_i|_{M^o} \equiv 0, \, 1 \leq i \leq k$, and by density
$\varphi_i \equiv 0, 1 \leq i \leq k.$ \qed

\begin{cor} 
The restrictions of $d\ell_1, d\ell_2, ... , d\ell_{n-3}$ to $M
\subset M_r$ are independent over $C^\8(M)$.
\end{cor}

\begin{rem}
Since \ $\ell$ \ is \ onto, \ if \ there \ exists \ $\Phi \in
C^{\8}(B_r)$ \ such \ that \ $\Phi(\ell_1(x),...,\ell_k(x)) \equiv 0$,
then $\Phi \equiv 0$.
\end{rem}

We conclude this chapter by giving a second proof that the bending flows on
disjoint diagonals commute.  Since $M_r^o$ is dense in $M_r$, it suffices to
prove

\begin{lem}
$\Psi_i^s(\Psi_j^t(P))=\Psi_j^t(\Psi_i^s(P))$, for $P\in M_r^o$.
\end{lem}

\proof  We assume $i>j$.  We observe that the diagonals $d_{1i}$ and $d_{1j}$
divide the surface $S$ into three polyhedral ``flaps'', $I, \;II, \;III$ (the
boundary of $I$ contains $e_1$, the boundary of $II$ contains $e_i$, and the
boundary of $III$ contains $e_j$).  Let $R_i^s$ and $R_j^t$ be the one
parameter groups of rotations around $d_{1i}$ and $d_{1j}$, respectively.  We
first record what $\Psi_i^s \circ \Psi_j^t$ does to the flaps.
\begin{align*}
&\Psi_i^s \circ \Psi_j^t(I)  =  R_i^s R_j^t(I) \\
&\Psi_i^s \circ \Psi_j^t(II)  = R_i^s(II)  \\
&\Psi_i^s \circ \Psi_j^t(III)  = III
\end{align*}

Now we compute what $\Psi_j^t \circ \Psi_i^s$ does to the flaps.  The point is,
after the bending on $d_{1i}$, the diagonal $d_{1j}$ moves $R_i^s d_{1j}$.
Hence, the next bending rotates $I$ around $R_i^s d_{1j}$.  Hence, the next
bending curve is $R_i^s \circ R_j^t \circ R_i^{-s}$.  We obtain
\begin{align*}
&\Psi_j^t \circ \Psi_i^s(I) = (R_i^s R_j^t R_i^{-s}) R_i^s(I) \\
&\Psi_j^t \circ \Psi_i^s(II) = R_i^s (II) \\
&\Psi_j^t \circ \Psi_i^s(III) = III.
\end{align*}
\qed

\section{Symplectomorphism of $M_r(\E^3)$ and $M_r(\h^3)$}

Recall $r$ is not on a wall of $D_n$.  Then by Theorem \ref{KKN} of this paper,
the
hyperbolic Gauss map $\ga=\ga_h:M_r(\h^3) \to Q_{sst}(r)$ is a diffeomorphism.
Moreover by Theorem 2.3 of \cite{KM2}, the Euclidean Gauss map $\ga_e:M_r(\E^3)
\to Q_{sst}(r)$ is also a diffeomorphism.  We obtain

\begin{thm}
Suppose $r$ is not on a wall of $D_r$, then the composition $\ga_h^{-1} \circ
\ga_e:M_r(\E^3) \to M_r(\h^3)$ is a diffeomorphism.
\end{thm}

\begin{rem}
The result that $M_r(\E^3)$ and $M_r(\h^3)$ are (noncanonically) diffeomorphic
was obtained by
\cite{Sa}.
\end{rem}

It does not appear to be true that $\ga_h^{-1} \circ \ga_e$ is a
symplectomorphism.

\subsection{A Formula of Lu} \label{form}

In the next several sections we will prove that $M_r(\h^3)$ is
symplectomorphic to $M_r(\E^3)$.

We first define a family of nondegenerate Poisson structures $\pi_\eps$,
$\eps\in[0,1]$, on the 2-sphere, $S^2 \simeq K/T$.  Letting $\om_\eps$ be the
corresponding family of symplectic forms we show the cohomology
classes $[\om_\eps]$ of $\om_\eps$ in $H^2(S^2)$ are constant.

Fix $\la \in \R_+$ and $\La = X \wedge Y \in \wedge^2 \k$, where $X=\half \left(
\begin{smallmatrix}
0 & 1 \\ -1 & 0 \end{smallmatrix} \right) $ and $Y= \half \left(
\begin{smallmatrix} 0 & i \\ i & 0 \end{smallmatrix} \right)$. The following
family of Poisson structures $\pi_\eps$ on
$K/T \simeq S^2$ for $\eps\in (0,1]$ are due to J.-H. Lu \cite{Lu2}.
$$
\pi_{\eps}=\eps[\pi_\8 - \tau (\eps) \pi_0]
$$
where $\pi_\8 = p_*\pi_K= p_*(dL_k \La - dR_k \La)$, $\tau(\eps) = \frac{1}{1-e^{4\eps\la}}$, and
$\pi_0 = 2 \, dL_k \, \La$.  Here $p:K \to K/T$ is the projection map. Then
$$
\pi_\eps (k) = \eps(dL_k \La - dR_k \La) -\frac{2\eps}{1- e^{4\eps\la}}\, dL_k \La.
$$

\begin{lem}
$$
\lim_{\eps\to0} \pi_\eps = \frac{1}{4\la}\pi_0
$$
\end{lem}
\proof  The proof of the lemma is a simple application of L'H\^{o}pital's rule.

\begin{lem} \label{abc}
$\pi_\eps$ is nondegenerate for $\eps\in [0,1]$.
\end{lem}

We will prove Lemma \ref{abc} in Proposition \ref{Symp}, where we show
$(K/T,\pi_\eps)$ is symplectomorphic to a
symplectic leaf of the Poisson Lie group $(B_\eps, \ha{\pi}_{B_\eps})$

We leave it to the reader to verify the Poisson structures on
$S^2$ can be written $$ \pi_\8 = \frac{1}{2}(1+\al^2+\be^2)
\frac{\D}{\D \al} \wedge \frac{\D}{\D \be}. $$ and $$ \pi_0 =
\frac{1}{2}(1+\al^2+\be^2)^2 \frac{\D}{\D \al} \wedge \frac{\D}{\D
\be} $$ where $(\al,\be)$ are coordinates obtained by
stereographic projection with respect to the north pole (see
\cite{LW}).  $\pi_\8$ is the Bruhat-Poisson structure on $K/T$. We
now let $\om_\eps$ be the symplectic form obtained by inverting
$\pi_\eps$ (this is possible since $\pi_\eps$ is nondegenerate).

$$ 
\om_\eps = \frac{- d\al \wedge d\be}{\eps(\half (1+\al^2+\be^2)
- \half \tau(\eps)(1+\al^2+\be^2)^2)} \;\;, \eps\in (0,1] 
$$ 
Let $\om_0$ be the limiting symplectic structure $$ \om_0 = -8\la
\frac{d\al\wedge d\be}{(1 + \al^2 + \be^2)^2} $$

\begin{lem}
$$
\int_{\R^2} \om_0 = -8\pi\la
$$
\end{lem}
\proof
\begin{eqnarray*}
\int_{\R^2} \om_0 & = & \int_{\R^2} -8\la \frac{d\al\wedge d\be}{(1 + \al^2 +
    \be^2)^2} \\
    & = & -8 \la \int_{\th=0}^{\th=2\pi} \int_{r=0}^{r=\8} \frac{ r\; dr\wedge
    d\th}{(1+ r^2)^2} \\
    & = & -16 \pi \la \int_{u=1}^{u=\8} \frac{(1/2) du}{u^2} \\
    & = & -8 \pi \la
\end{eqnarray*}
\qed

\begin{lem}
$$
\int_{\R^2} \om_\eps = -8\pi\la, \:\eps\in(0,1]
$$
\end{lem}
\proof  Note that $\tau(\eps)<0$.   Then
\begin{eqnarray*}
    \int_{\R^2} \om_\eps & = & -\frac{2}{\eps}\int_{\R^2} \frac{d\al \wedge
d\be}{(1+\al^2+\be^2) - \tau(\eps)(1+\al^2+\be^2)^2} \\
    & = & -\frac{2}{\eps} \int_{\th=0}^{\th=2\pi} \int_{r=0}^{r=\8} \frac{r\;dr\wedge
d\th}{(1+r^2) - \tau(\eps)(1+r^2)^2} \\
    & = & -\frac{4\pi}{\eps} \int_{u=1}^{u=\8} \frac{(1/2)du}{u -
    \tau(\eps)u^2)} \\
    & = & -\frac{2\pi}{\eps}\; \log
    {\left|\frac{\tau(\eps)-1}{\tau(\eps)}\right|} \\
    & = & -\frac{2\pi}{\eps}\; \log(e^{4\eps\la}) \\
    & = & -8\pi \la
\end{eqnarray*}
\qed

\no We have proved the following

\begin{lem} \label{luref}
The cohomology classes $[\om_\eps]$ of $\om_\eps$  in $H^2(S^2)$ are constant.
\end{lem}

\begin{rem}
The previous lemma is a special case of Lemma 5.1 of \cite{GW}.
\end{rem}

\subsection{Symplectomorphism of ($\Sigma_\la(\eps),\pi_{B_\eps}$) and
($K/T,\pi_{\la,\eps}$)} \label{15}

In this section we obtain the Poisson structure $\pi_\eps$ from a deformed Manin
triple $(\g_\eps, \k, \b_\eps)$.

For $\eps >0$, we define the isomorphism $f_\eps:\g_\eps \to \g$ by $f_\eps =
\rho_\k + \eps\rho_\b$,
so that $f_\eps(X+\xi)=X+\eps\xi$ for $X\in \k$ and $\xi \in \b$.  We will
define a Lie bracket on $\g_\eps$ by the pullback of the Lie bracket on $\g$,
$[u,v]_\eps=f_{1/ \eps}[f_\eps u,f_\eps v]$.
We also define $\<,\>_\eps$ as the pullback of $\<,\>$.  Here $[,]$ and $\<,\>$
are the usual structures on $\g$.  We define $\B_\eps:\g \to \g^*$ as the map
induced by $\<,\>_\eps$.  To simplify notation, the subscripts will be dropped
when $\eps = 1$.

The following lemma gives us a formula for the Lie bracket on $\g_\eps$.
\begin{lem}
$[X+\al, Y+\be]_\eps = [X,Y] + \eps\rho_\k[X,\be] + \eps\rho_\k
[\al,Y] + \rho_\b[X,\be] + \rho_\b[\al,Y] + \eps[\al,\be]$, where $X,Y \in \k$
and $\al,\be \in \b$.
\end{lem}
\proof
\begin{eqnarray*}
[X+\al, Y+\be]_\eps & = & f_{1/ \eps}[f_\eps (X+\al), f_\eps (Y+\be)]\\
    & = & f_{1/ \eps} [X+ \eps\al, Y+\eps\be]\\
    & = & f_{1/ \eps} \{[X,Y] + \eps[X,\be] + \eps[\al,Y] +
\eps^2[\al,\be]\}\\
    & = & f_{1/ \eps} \{[X,Y] \!+\! \eps\rho_\k[X,\be] \!+\! 
\eps\rho_\b[X,\be] \!+ \!\eps\rho_\k[\al,Y]\! +\! \eps\rho_\b[\al,Y] \! +
\!    \eps^2[\al,\be]\}\\
    & = & [X,Y] + \eps\rho_\k[X,\be] + \eps\rho_\k[\al,Y] + \rho_\b[X,\be] +
\rho_\b[\al,Y] + \eps[\al,\be]
\end{eqnarray*}
\qed

\no We leave it to the reader to check
\begin{lem}
$\<,\>_\eps = \eps\<,\>$
\end{lem}

Let $G_\eps$ be the simply-connected Lie group with Lie algebra
$\g_\eps$.  Let $F_\eps : G_\eps \to G$ be the isomorphism induced
by $f_\eps$. We have a commutative diagram of isomorphisms.

$$
  \begin{CD}
    \g_\eps @>f_\eps>> \g \\
    @V{\exp^\eps}VV  @VV{\exp}V \\
    G_\eps @>F_\eps>> G
  \end{CD}
$$

Let $x\in\g_\eps$. We use the identity map to identify $\g$ and
$\g_\eps$ as vector spaces. In what follows, we will make frequent use
of
\begin{lem}
$\ha{Ad}(\exp^\eps x) = Ad (\exp{\eps x})$ as elements in $GL(\b)$ for all $x
\in \b_\eps = \b$.  Here $\ha{Ad}$ denotes the adjoint action of
$G_\eps$ on $\g_\eps$.
\end{lem}
\proof By \cite[pg. 114]{Wa},
\begin{eqnarray*}
\ha{Ad}(\exp^\eps x) & = & e^{\ha{ad}x} \\
	& = & e^{\eps{ad}x} \\
	& = & e^{ad(\eps x)} \\
	& = & Ad (\exp {\eps x})
\end{eqnarray*}
\qed

Given our deformed Manin triple on $G_\eps$, $(\g_\eps, \k, \b_\eps)$, we will
construct a
Poisson structure $\ha{\pi}_{B_\eps}$ on $B_\eps$, the simply-connected Lie
group with Lie algebra, $\b_\eps$.  We will denote all quantities associated to
the deformed Manin triple with a hat \:$\ha{     }$ \;.

We define the Poisson Lie structure on $B_\eps$ by the Lu-Weinstein Poisson
tensor \cite{LW}
$$
\ha{\pi}_{B_\eps}(b)(\ha{dR}^*_{b^{-1}} \al_X,\ha{dR}^*_{b^{-1}}
\al_Y)=\<\rho_\k(\ha{Ad}_{b^{-1}}\B_\eps^{-1}(\al_X)),\rho_\b(\ha{Ad}_{b^{-1}}\B_
\eps^{-1} (\al_Y))\>_\eps
$$
where $\al_X,\al_Y\in \b_\eps^*$, $\al_X = \<X,\cdot\>_1$ and $\al_Y
= \<Y, \cdot\>_1$. 

\begin{rem}
Since $\lim_{\eps \to 0} \<,\>_\eps = 0$, it appears as if the limiting Poisson
structure $\lim_{\eps\to 0} \ha{\pi}_{B_\eps}$ will vanish.  However, we will
see in Proposition \ref{lim} that the limiting Poisson structure is associated
to the Manin triple $(\g_0, \k, \b_0)$ and $\<,\>_0 = \frac{d}{d\eps}\big|_{\eps=0} \<,\>_\eps$.
\end{rem}

We denote by $\pi_{B_\eps}$ the Poisson structure on $B_\eps$ using
the scaled bilinear form $\frac{1}{\eps}\<,\>_\eps=\<,\>$.  Then
$$
\pi_{B_\eps} (b)(\ha{dR}^*_{b^{-1}}\al_X,\ha{dR}^*_{b^{-1}}
\al_Y)=\<\rho_\k(\ha{Ad}_{b^{-1}}\B_1^{-1}(\al_X)),
\rho_\b(\ha{Ad}_{b^{-1}}\B_1^{-1}(\al_Y))\>
$$
where $\al_X,\al_Y\in \b_\eps^*$.  For the following we will let
$X_\eps=\B^{-1}_\eps (\al_X)\in \k$, again dropping the subscript when $\eps = 1$,
so that $X_\eps = \frac{1}{\eps}X$.

\begin{lem} \label{piBe}
$\ha\pi_{B_\eps}=\frac{1}{\eps}\pi_{B_\eps}$
\end{lem}
\proof
\begin{eqnarray*}
\ha{\pi}_{B_\eps}(b)(\ha{dR}^*_{b^{-1}}\al_X,\ha{dR}^*_{b^{-1}} \al_Y) & = &
\<\rho_\k(\ha{Ad}_{b^{-1}}X_\eps),\rho_\b(\ha{Ad}_{b^{-1}}Y_\eps)\>_\eps \\
    & = & \<\rho_\k(\ha{Ad}_{b^{-1}}\frac{1}{\eps}X),\rho_\b(\ha{Ad}_{b^{-1}}
    \frac{1}{\eps}Y)\>_\eps \\
     & = & \frac{1}{\eps}\<\rho_\k(\ha{Ad}_{b^{-1}}X),\rho_\b(\ha{Ad}_{b^{-1
     }}Y)\>
\end{eqnarray*}
\qed

\begin{prop} \label{lim}
$\lim_{\eps\to0}\ha{\pi}_{B_\eps}(b)(\ha{dR}^*_{b^{-1}} \al_X,\ha{dR}^*_{b^{-1}} \al_Y)=-\<\log
b, [X,Y]\>$.
\end{prop}
\proof
\begin{eqnarray*}
\le\ha{\pi}_{B_\eps}(b)(\ha{dR}^*_{b^{-1}} \al_X,\ha{dR}^*_{b^{-1}}
\al_Y) & = & \le\frac{1}{\eps}
    \<\rho_\k(\ha{Ad}_{b^{-1}}X),\rho_\b(\ha{Ad}_{b^{-1}}Y)\> \\
    & = & \le \frac{1}{\eps}2Im\, tr(\rho_\k
    (Ad_{e^{-\eps \log b}}X) \rho_\b(Ad_{e^{-\eps \log b}}Y)) \\
    & = & \le 2Im \,tr(\rho_\k
    (Ad_{e^{-\eps \log b}}X) \rho_\b(-\log b \, Y +Y \log b)) \\
    & = & - 2Im \,tr(X \rho_\b[\log b, Y ]) \\
    & = & - 2Im \,tr(X [\log b, Y ]) \\
    & = & - 2Im \,tr(\log b [X, Y ]) \\
    & = & - \<\log b, [X, Y ]\>
\end{eqnarray*}
\qed

\begin{rem}
Before stating the next corollary, note that the limit Lie algebra
$\b_0$ is abelian whence the limit Lie group $B_0$ is abelian.  Hence,
$\exp_0:\b_0 = T_o(B_0) \to B_0$ is the canonical identification of
the vector space $B_0$ with its tangent space at the origin.  Hence,
$\exp_0$ is an isomorphism of Lie groups and $\exp_0^*$ carries
invariant 1-forms on $B_0$ to invariant 1-forms on $\b_0$.
\end{rem}

\begin{cor}
$\lim_{\eps \to 0} \ha{\pi}_{B_\eps}$ is the negative of the Lie
Poisson structure on $\k^* \simeq \b_0$ transferred to $B_0$ using the
exponential map on the vector space $B_0$.
\end{cor}
\proof The proof is left to the reader.

\begin{rem}
$(K, \eps \pi_K)$ is the dual Poisson Lie group of $(B_\eps,
\ha{\pi}_{B_\eps})$.
\end{rem}

We will denote the dressing action of $K$ on $B_\eps$ by
$\ha{D}_\eps^\ell$ and the infinitesimal dressing action of $\k$ on $B_\eps$ by
$\ha{d_\eps^\ell}$.  By definition $\ha{d_\eps^\ell}(b)(X) =
\ha{\pi}_{B_\eps}(\cdot,\al_X)$.  We then have the following.

\begin{lem}
$\ha{d_\eps^\ell}(b)(X)=\frac{1}{\eps}d^\ell(b)(X)$
\end{lem}

\proof
Follows immediately from Lemma \ref{piBe} and the definition of
dressing action. 
\qed

\begin{rem}
$\le\ha{d_\eps^\ell}(b)(X)=ad^*(X)(\log b)$
\end{rem}

For the remainder of the section, fix $\la \in \R_+$ and $a=
\exp^\eps{\la H}\in B_\eps$, where $H = diag(1, -1) \in \mathfrak{a}_\eps$.
Let $\varphi_\eps:K \to \Sigma_\la^\eps \subset B_\eps$ be the map defined
by $\varphi_\eps(k)= \ha{D}_\eps^\ell(k)(a) = \rho_{B_\eps}(k*a)$, where
$\Sigma_\la^\eps$ is the
symplectic leaf through the point $a \in B_\eps$.  The map $\varphi_\eps$
induces a diffeomorphism from $K/T$ onto $\Sigma_\la^\eps$ which we will also
denote by $\varphi_\eps$.
\smallskip
\no Recall the family of Poisson tensors on $K/T$ given in \S \ref{form}
$$
\pi_{\la,\eps} = \eps(\pi_\8 - \tau(\eps\la) \pi_0).
$$

\begin{lem}
The map $\varphi_\eps : K/T \to B_\eps$ is $K$-equivariant, where
$(K,\eps\pi_K)$ acts on $(K/T, \pi_{\la,\eps})$ by left multiplication and
$B_\eps$ by the dressing action.
\end{lem}
\proof  $\varphi_\eps(g\cdot k) = \ha{D}_\eps^\ell(gk)(a) =
\ha{D}_\eps^\ell(g) (\ha{D}_\eps^\ell(k)(a)) =
g\cdot\varphi_\eps(k)$. \qed

\begin{rem}
The action of $(K,\eps\pi_K)$ on $(K/T, \pi_{\la,\eps})$ by left
multiplication is a Poisson action.
\end{rem}

Since $K/T$ is a symplectic manifold, there is a momentum map for
the action of $K$ on $K/T$, see \cite[Theorem 3.16]{Lu}.  We will see
as a consequence of Proposition \ref{Symp}

\begin{lem}
The momentum map for the action of $(K,\eps\pi_K)$ on $(K/T,\pi_{\la,\eps})$ is
$\varphi_\eps$.
\end{lem}

\begin{prop} \label{Symp}
The map $\varphi_\eps$ induces a symplectomorphism from $(K/T,
\pi_{\la,\eps})$ to $\: \: \: $ $(\Sigma^\eps_\la, \ha{\pi}_{B_\eps})$.
\end{prop}

\proof   Since the $K$-actions on $K/T$ and $\Sigma^\eps_\la$ are Poisson and
the map $\varphi_\eps : K/T \to \Sigma^\eps_\la$ is a $K$-equivariant
diffeomorphism, if $(d\varphi_\eps)_e(\pi_{\la,\eps}(e)) =
\ha\pi_{B_\eps}(a)$ then it follows that
$(d\varphi_\eps)_k(\pi_{\la,\eps}(k)) = 
\ha\pi_{B_\eps}(\varphi(k))$ for
all $k \in K/T$.

We will need the following lemmas to prove the proposition.
We let $E= \left(
\begin{smallmatrix}
0 & 1 \\
0 & 0
\end{smallmatrix}
\right)\in \b$ and $\Lambda = E\wedge iE \in \b \wedge \b$.  If we set
$\ha{\pi}_\La(b)
= \frac{1}{\eps}(\ha{dL}_b \,\La - \ha{dR}_b \,\La)$, we then have the
following.

\begin{lem}
$\ha{\pi}_{B_\eps}|_a = \frac{1}{2} \ha{\pi}_\La|_a$ for $a=\exp^\eps {\la H}$.
\end{lem}
\proof
Let $X = \left( \begin{smallmatrix} si & u \\ -\bar{u} & -si \end{smallmatrix}
\right)$,  $Y = \left( \begin{smallmatrix} ti & v \\ -\bar{v} & -ti
\end{smallmatrix} \right) \in \k$
\begin{eqnarray*}
\ha{\pi}_{B_\eps}(a)(\ha{dR}^*_{a^{-1}} \al_X,\ha{dR}^*_{a^{-1}} \al_Y) & = & \frac{1}{\eps}
    \<\rho_\k(\ha{Ad}_{a^{-1}}X),\rho_\b(\ha{Ad}_{a^{-1}}Y)\> \\
    & = & \frac{1}{\eps}
    \<X,\ha{Ad}_{a}\rho_\b(\ha{Ad}_{a^{-1}}Y)\>
\end{eqnarray*}
We can see,
\begin{eqnarray*}
\ha{Ad}_{a}\rho_\b(\ha{Ad}_{a^{-1}}Y) & = & \ha{Ad}_{a}\rho_\b\left(\ha{Ad}_{a^{-1}}
\left( \begin{smallmatrix} ti & v \\ -\bar{v} & -ti \end{smallmatrix} \right)\right)
\\
& = & \ha{Ad}_{a}\rho_\b
\left( \begin{smallmatrix} ti & e^{-2\eps\la}v \\ -e^{2\eps\la}\bar{v} & -ti
\end{smallmatrix} \right)
\\
& = & \ha{Ad}_{a}
\left( \begin{smallmatrix} 0 & (e^{-2\eps\la}-e^{2\eps\la})v \\ 0 & 0
\end{smallmatrix}
\right)
\\
& = & (1- e^{4\eps\la})
\left( \begin{smallmatrix} 0 & v \\ 0 & 0 \end{smallmatrix} \right)
\end{eqnarray*}
so that
\begin{eqnarray*}
\frac{1}{\eps}\<X,\ha{Ad}_{a}\rho_\b(\ha{Ad}_{a^{-1}}Y)\> & = & \frac{2}{\eps}
Im tr\left[\left( \begin{smallmatrix} si & u \\ -\bar{u} & -si \end{smallmatrix}
\right)
\left( \begin{smallmatrix} 0 & (1- e^{4\eps\la})v \\ 0 & 0 \end{smallmatrix}
\right)\right]
\\
    & = & -\frac{2}{\eps}(1- e^{4\eps\la})Im (\bar{u}v)\\
    & = & \frac{2}{\eps}(e^{4\eps\la} - 1)Im (\bar{u}v).
\end{eqnarray*}

\no If we evaluate the right-hand side of the above formula we see

\begin{eqnarray*}
\half \ha{\pi}_\La(a)(\ha{dR}^*_{a^{-1}} \al_X,\ha{dR}^*_{a^{-1}}
\al_Y) & = & \frac{1}{2\eps} [\al_X \wedge \al_Y (\ha{Ad}_a E,
\ha{Ad}_a iE) - \al_X \wedge \al_Y(E, iE) \\   
    & = & \frac{2}{\eps}[e^{4\eps\la} Im(\bar{u}v) - Im(\bar{u}v)] \\
    & = & \frac{2}{\eps}(e^{4\eps\la} -1) Im(\bar{u}v) \\
    & = & \ha{\pi}_{B_\eps}(a)(\ha{dR}^*_{a^{-1}}
\al_X,\ha{dR}^*_{a^{-1}} \al_Y) 
    \end{eqnarray*}
\qed

\no We then have the following.

\begin{cor}
$\ha\pi_{B_\eps}(a) = \frac{1}{2\eps}(1 - e^{-4\eps\la}) \ha{dL}_a \, (E \wedge
iE)$
\end{cor}
\proof
\begin{eqnarray*}
\ha\pi_{B_\eps}(a) & = & \frac{1}{2\eps}[\ha{dL}_a (E \wedge iE) - \ha{dR}_a (E \wedge
iE) ] \\
    & = & \frac{1}{2\eps}\ha{dL}_a[E \wedge iE - \ha{Ad}_{a^{-1}}\, (E \wedge iE)] \\
    & = & \frac{1}{2\eps}\ha{dL}_a[E \wedge iE -  e^{-4\eps\la} (\,E \wedge iE)] \\
    & = & \frac{1}{2\eps}(1 - e^{-4\eps\la}) \ha{dL}_a \, (E \wedge iE).
\end{eqnarray*}
\qed

The diffeomorphism $ \varphi_\eps:K/T \to \Sigma_\la^\eps$ gives us
$(d\varphi_\eps)_e: \k / \mathfrak{t}
\to T_a\Sigma_\la^\eps \subset T_a B_\eps$ defined by
$(d\varphi_\eps)_e(\xi)=\frac{1}{\eps}
\ha{dL}_a \rho_\b(\ha{Ad}_{a^{-1}}\xi)$.  Now let $X=\half \left(
\begin{smallmatrix}
0 & 1 \\ -1 & 0 \end{smallmatrix} \right) $ and $Y= \half \left(
\begin{smallmatrix} 0 & i \\ i & 0 \end{smallmatrix} \right)$ as in \S
\ref{form}, then
$$
(d\varphi_\eps)_e(X)=\frac{1}{2\eps}(e^{-2\eps\la}-e^{2\eps\la})\ha{dL}_a\,
E\: \: \: 
\textnormal{and}
\: \: \:(d\varphi_\eps)_e(Y)=\frac{1}{2\eps}(e^{-2\eps\la}-
e^{2\eps\la})\ha{dL}_a\, iE.
$$
It then follows that
\begin{lem}
$(d\varphi_\eps)_e(\pi_{\la,\eps}(e)) = \ha\pi_{B_\eps}(a)$
\end{lem}
\proof
\begin{eqnarray*}
(d\varphi_\eps)_e(\pi_{\la,\eps}(e)) & = & (d\varphi_\eps)_e(\eps(\pi_\8(e) -
\tau(\eps\la)\pi_0(e))) \\
    & = & \eps(d\varphi_\eps)_e(\pi_\8(e)) - \eps\tau(\eps\la)
    (d\varphi_\eps)_e(\pi_0(e)) \\
    & = & 0 - 2\eps\tau(\eps\la) (d\varphi_\eps)_e(X\wedge Y) \\
    & = & -\frac{1}{2\eps}\tau(\eps\la) (e^{-2\eps\la}-e^{2\eps\la})^2 \ha{dL}_a
    (E\wedge iE) \\
    & = & -\frac{1}{2\eps}(e^{-4\eps\la} - 1 ) \ha{dL}_a (E \wedge iE) \\
    & = & \ha\pi_{B_\eps}(a).
\end{eqnarray*}
\qed 

\smallskip \no This completes the proof of Proposition \ref{Symp}.
\qed

\medskip

We can next look at the product $(K/T)^n$.  We give $(K/T)^n$ the
product Poisson structure $\pi_{\la,\eps} = \pi_{\la_1,\eps} +
\cdots + \pi_{\la_n,\eps}$.  Define the map $$ \t{\Phi}_\eps:
(K/T)^n \to \Sigma_{\la_1}^\eps \times \cdots \times
\Sigma_{\la_n}^\eps $$ given by $$ \t{\Phi}_\eps(k_1,...,k_n) =
(\varphi^\eps_1(k_1), ... , \varphi^\eps_n(k_n)) =
(\ha{D}^\ell_\eps(k_1)(a_{\la_1}), ... ,
\ha{D}^\ell_\eps(k_n)(a_{\la_n})) $$ where $\varphi^\eps_i(k_i) =
\ha{D}^\ell_\eps(k_i)(a_{\la_i})$ and  $a_{\la_i} = \exp^\eps
(\la_i H) \in B_\eps$.  We note the map $\t{\Phi}_\eps: (K/T)^n \to
\Sigma_{\la_1}^\eps \times \cdots \times \Sigma_{\la_n}^\eps$ is a
symplectomorphism.

We leave the proof of the following lemma to the reader.
\begin{lem} \label{action}
The action of $K$ on $(K/T)^n$ given by
$$
k \circ (k_1, ... , k_n) = (kk_1, \rho_K(k \varphi^\eps_1(k_1))k_2,...,
\rho_K(k \varphi^\eps_1(k_1)\cdots \varphi^\eps_{n-1}(k_{n-1}))k_n)
$$
is the pull back under $\t{\Phi}_\eps$ of the $\eps$-dressing action on
$\Sigma_{\la_1}^\eps \times \cdots \times \Sigma_{\la_n}^\eps \subset B_\eps^n$.
\end{lem}

The momentum map for the action of $(K,\eps\pi_K)$ on
$((K/T)^n,\pi_{\la,\eps})$ is $$ \t\Psi_\eps: (K/T)^n \to B_\eps
$$ where $$ \t\Psi_\eps(k_1, .. , k_n) = \varphi^\eps_1(k_1) *
\cdots
* \varphi^\eps_{n-1}(k_{n-1}). $$

\subsection{The $\eps$-dressing orbits are small spheres in hyperbolic 3-space}

Let $b$ be the Killing form on $\g$ divided by 8.  We have normalized $b$ so
that the
induced Riemannian metric $(,)$ on $G/K$ has constant curvature -1.  We let
$b_\eps =
f_\eps^*b$, hence $b_\eps$ is the Killing form on $\g_\eps$.  Then $(,)_\eps =
F_\eps^*(,)$ is
the induced Riemannian metric on $G_\eps/K$ and $G_\eps/K$ has constant
curvature -1 (since
$F_\eps$ is an isometry).  We will call $(,)_\eps$ the hyperbolic metric on
$G_\eps/K$.

The map $\zeta:B_\eps \to G_\eps/K$ given by $\zeta(b)=b*K$ is a diffeomorphism
that
intertwines the $\eps$-dressing orbits of $K$ on $B_\eps$ with the natural $K$
action on
$G_\eps/K$ given by left multiplication (using the multiplication in $G_\eps$).
We
abbreviate the identity coset $K$ in $G_\eps/K$ to $x_0$ and use the same letter
for the
corresponding point in $G/K$.  We have

\begin{lem} \label{orbit}
The image of the $\eps$-dressing orbit $\Sigma_\la^\eps$ under $\zeta$ is the
sphere around
$x_0$ of radius $\eps\la$.
\end{lem}
\proof  Let $d_\eps$ be the Riemannian distance function on
$G_\eps/K$ and $d$ the Riemannian distance function on $G/K$.  We
have
\begin{eqnarray*}
d_\eps(x_0, \exp_{x_0}^\eps \la H) & = & d(x_0, F_\eps \exp_{x_0}^\eps \la H) \\
    & = & d(x_0, \exp_{x_0} f_\eps(\la H)) \\
    & = & d(x_0, \exp_{x_0} \eps \la H) \\
    & = & \eps\la
\end{eqnarray*}
\qed

\subsection{The family of symplectic quotients} \label{17}
In this section we will continue to use the notation of
\S\ref{15}.  Let $p: E = (K/T)^n \times I \to I$ be a projection.
Here we define $I=[0,1]$.
We let $T^{vert}(E) \subset T(E)$ be the tangent space to the
fibers of $p$. Hence $\bigwedge^2 T^{vert}(E)$ is a subbundle of
$\bigwedge^2 T(E)$. We define a Poisson bivector $\pi$ on E by
$\pi(u,\eps) = \pi_{\la,\eps} |_u$.  $\pi$ is a section of
$\bigwedge^2 T^{vert}(E)$.

Let $S \subset B^n \times I$ be defined by $S= \{(b,\eps) | b\in \Sigma^\eps_\la
\}$.  We let $(K, \eps\pi_K)$ act on $(K/T)^n \times I$ by $k\cdot(u,\eps)=
(k\circ u,\eps)$, where $\circ$ is the action given in Lemma \ref{action},
and act on $S$ by $k\cdot(b,\eps)=(\ha{D}^\ell_\eps(k)(b), \eps)$.
We then define the map $\Phi : E \to S$ by $\Phi(u,\eps) =
(\t{\Phi}_\eps(u),\eps)$ which  is a K-equivariant diffeomorphism.  We also
define $\Psi : E \to B$ by $\Psi(u,\eps) = \t{\Psi}_\eps(u)$.

\begin{rem}
$\Psi|_{p^{-1}(\eps)}$ is the momentum map for the Poisson action of
$(K,\eps \pi_K)$ on $(K/T)^n \times \{\eps \}$.
\end{rem}

We need some notation.  Suppose $n \geq m$, $F : \R^{n+1} \to
\R^m$ is a smooth map, and $0\in \R^m$ is a regular value of $F$.
Let $M = F^{-1}(0)$.  Write $\R^{n+1} = \R^n \times \R$ with $x
\in \R^n, \; t\in \R$.  Let $p: \R^{n+1} \to \R$ be the projection
onto the $t$-line.  The next lemma is taken from \cite{Sa}.

\begin{lem}
Let $(x,t)\in M$.  Suppose $\frac{\partial F}{\partial x}\big|_{(x,t)}$ has
maximal rank $m$.  Then $dp|_{(x,t)}:T_{(x,t)}(M) \to T_t(\R)$ is onto.
\end{lem}
\proof  It suffices to construct a tangent vector $v\in T_{(x,t)}(\R^{n+1})$
satisfying

\begin{itemize}
\item[(i)]  $v \in ker\; dF|_{(x,t)}$
\item[(ii)]  $v = \sum_{i=1}^n c_i \frac{\partial}{\partial x_i} +
\frac{\partial}{\partial t}$
\end{itemize}

Put $c = (c_1, ... , c_n)$ and write the Jacobian matrix $dF|_{(x,t)}$ as
$(A,b)$ where $A$ is the $m$ by $n$ matrix given by $A = \frac{\partial
F}{\partial x}|_{(x,t)}$ and $b$ is the column vector of length $m$ given by $b=
\frac{\partial F}{\partial t}|_{(x,t)}$.  We are done if we can solve
$$
Ac+b=0.
$$
But since $A: \R^n \to \R^m$ is onto we can solve this equation.
\qed

\begin{rem}
We need to generalize to the case in which $\R^{n+1}$ is replaced
by the closed half-space $\bar{H} = \{(x,t):\; x\in\R^n, \; t \geq
0\}$ and the $t$-line by the closed half-line.  Given $(x,0) \in
\D M$ we wish to find $v= \sum_{i=1}^n c_i \frac{\D}{\D x_i} +
\frac{\D}{\D t}$ with $dF|_{(x,o)} (v) = 0$ (so $v$ is in the
tangent half-space to $M$ at $(x,0) \in \D M$). The argument is
analogous to that of the lemma and is left to the reader.
\end{rem}

\begin{cor}
Suppose $M$ is compact and for all $(x,t)$ and further that $\frac{\partial
F}{\partial x}|_{(x,t)}$ has maximal rank for all $(x,t) \in M$.  Then $p:M\to
\R$ is a trivial fiber bundle.
\end{cor}
\proof  $p$ is proper since $M$ is compact.  This is the Ehresmann fibration
theorem \cite[8.12]{BJ}.
\qed

\begin{rem}
We leave to the reader the task of extending the corollary to the case where
the $t$-line is replaced by the closed $t$ half-line.
\end{rem}

We now return to our map $\Psi : (K/T)^n \times I \to B$.  We have
$$
\Psi(u,\eps) = \varphi_1^\eps(u_1) * \cdots * \varphi_n^\eps(u_1).
$$
Let $\eps > 0$.  We apply Lemma \ref{regval} to deduce that $1\in B$ is a
regular value for
$u \to \Psi(u,\eps)$ (recall we have assumed $r$ is not on a wall of $D_n$).
Now let $\eps = 0$.  It is immediate (see \cite{KM2}) that $0\in \R^3$ is a
regular value of $u \to \Psi(u,\eps)$ (again because $r$ is not on a wall of
$D_n$).  We obtain

\begin{lem}
$p: \Psi^{-1}(1) \to I$ is a trivial fiber bundle.
\end{lem}

Now let $\mathcal{M} = \Psi^{-1}(1)/K$.  We note that $p$ factors through the
free action of $K$ on $\Psi^{-1}(1)$ and  we obtain a fiber bundle
$\bar{p}:\mathcal{M} \to I$.  This gives the required family of symplectic
quotients.

\begin{prop}
$\bar{p} : \M \to I$ is a trivial fiber bundle.
\end{prop}

\begin{rem}
$\bar{p}^{-1}(0) = M_r(\E^3)$ and $\bar{p}^{-1}(1) = M_r(\h^3)$ and we may
identify $\M$ with the product $M_r(\E^3) \times I$.
\end{rem}

We now give a description of the symplectic form along the fibers of
$\bar{p}$.  Recall that if $\pi : E \to B$ is a smooth fiber bundle then the
relative forms on $E$ are the elements of the quotient of $A^\bullet
(E)$ by the 
ideal generated by elements of positive degree in $\pi^* A^\bullet(B)$.  Note the restriction of a relative
form to a fiber of $\pi$ is well-defined and the relative forms are a
differential graded-commutative algebra with product and differential induced by
those of $A^\bullet(E)$.

Lu's one parameter family of forms $\om_\eps$ of \S\ref{form} induces a relative
2-form $\om_\eps$ on $K/T \times I$ which is relatively closed.  By taking sums
we obtained a relative 2-form $\t{\om}_\eps$ on $(K/T)^n \times I$ and by
restriction and projection a relative 2-form $\bar{\om}_\eps$ on $\M$.  Clearly
$\bar{\om}_\eps$ is relatively closed and induces the symplectic form along the
fibers of $\bar{p}: \M \to I$.

We let $[\bar{\om}_\eps]$ be the class in $H^2(p^{-1}(\eps))$ determined by
$\bar{\om}_\eps$.

\subsection{$[\bar{\om}_\eps]$ is constant and Moser's Theorem}

To complete our proof we need to review the $i$-th cohomology bundle associated
to a smooth fiber bundle $\pi:E\to B$.  The total space $\mh^i_E$ of the $i$-th
cohomology bundle is given by $\mh_E^i = \{(b,z):\; b\in B, \; z\in
H^i(p^{-1}(b)) \}$.  We
note that a trivialization of $E|_U$ induces an isomorphism between
$\mh^i_{E|_U}$ and $\mh^i_{U\times F}$.  But $\mh^i_{U\times F} = \{(x,z) :\;
x\in U,\; z
\in H^i(F)\}$, whence $\mh^i_{U \times F} = U \times H^i(F)$.  It is then clear
that $\mh^i_E$ is a vector bundle over $B$ with typical fiber $H^i(F)$.  We next
observe that the action of the transition functions of $E$ on $H^i(F)$ induce
the transition functions of $\mh^i_E$.  Hence if we trivialize $E$ relative to a
covering $\U=\{U_i : i \in I\}$ such that all pairwise intersections are
contractible then the corresponding transition functions of $\mh^i_E$ are
constant.  Hence $\mh^i_E$ admits a flat connection called the Gauss-Manin
connection.  We observe that a cross-section of $\mh^i_E$ is parallel for the
Gauss-Manin connection if when  expressed locally as an  element  of
$\mh^i_{U\times F}$ as above it corresponds to a constant map from $U$ to
$H^i(F)$.

\begin{rem}
If $\tau$ is a relative $i$-form on $E$ which is relatively closed then it gives
rise to a cross-section $[\tau]$ of $\mh^i_E$ such that $[\tau](b)$ is the de
Rham cohomology class of $\tau(b)|_{\pi^{-1}(b)}$.
\end{rem}

We now consider the relative 2-form $\bar{\om}_\eps$ on
$\mathcal{M}$.  The form $\bar{\om}_\eps$ is obtained from the
corresponding form $\t{\om}_\eps$ on $(K/T)^n \times I$ by first
pulling $\t{\om}_\eps$ back to $\Psi^{-1}(1)$ then using the
invariance of $\t{\om}_\eps$ under $K$ to descend $\t{\om}_\eps$
to $\bar{\om}_\eps$. We observe that $[\t{\om}_\eps]$ (reps.
$[\bar{\om}_\eps]$) is a smooth section of $\mh^2_{(K/T)^n\times
I}$ (resp. $\mh^2_{\mathcal{M}}$)).

We obtain a diagram of second cohomology bundles with connection
$$
  \begin{CD}
  \mh^2_{(K/T)^n\times I} @>{i^*}>>  \mh^2_{\Psi^{-1}(1)} @< \pi^* <<
\mh^2_{\mathcal{M}} \\
  \end{CD}
$$
where $i: \Psi^{-1}(1) \to (K/T)^n\times I$ is the inclusion and
$\pi:\Psi^{-1}(1) \to \mathcal{M}$ is the quotient map.  We have
$$
\pi^*[\bar{\om}_\eps]= i^*[\t{\om}_\eps].
$$
\begin{prop}
$[\bar{\om}_\eps]$ is parallel for the Gauss-Manin connection on
$\mh^2_\mathcal{M}$.
\end{prop}

\proof  \ \ 
By Lemma \ref{luref}, $[\t{\om}_\eps]$ is \ parallel for the \ Gauss-Manin
connection on $(K/T)^n \times I$.
Hence $i^*[\t{\om}_\eps]$ is parallel for the Gauss-Manin connection on
$\mh^2_{\Psi^{-1}(1)}$.  But an elementary spectral sequence argument for the
bundle $K \to \Psi^{-1}(1) \to \mathcal{M}$ shows that $\pi^* :
\mh^2_{\mathcal{M}} \to \mh^2_{\Psi^{-1}(1)}$ is a bundle monomorphism.  Hence
if
$\pi^*[\bar{\om}_\eps]$ is parallel, so is $[\bar{\om}_\eps]$.
\qed

\begin{cor}
The cohomology class of $\bar{\om}_\eps$ is constant relative to any
trivialization of $p:\mathcal{M} \to I$.
\end{cor}

We now complete the proof of symplectomorphism by applying a version of Moser's Theorem
\cite{M} with $M = M_r(\E^3)$.  For the benefit of the reader we will
state and
prove the version of Moser's Theorem we need here.

\begin{thm}
Suppose $\om_\eps$ is a smooth one-parameter family of symplectic forms on a
compact smooth manifold $M$. Suppose the cohomology class $[\om_\eps]$ of
$\om_\eps$ in $H^2(M)$ is constant.  Then there is a smooth curve $\phi_\eps$ in
Diff($M$) with  $\phi_0 = id_M$ such that
$$
\om_\eps = \phi_\eps^* \om_0.
$$
\end{thm}
\proof  Choose a smooth one-parameter family of 1-forms $\tau_\eps$ such that
$$
\frac{d\om_\eps}{d\eps} = - d \tau_\eps.
$$
(We may choose $\tau_\eps$ {\em smoothly} by first choosing a Riemannian metric
then taking $\tau_\eps$ to be the coexact primitive of
$\frac{d\om_\eps}{d\eps}$ -
here we use the compactness of $M$).

Let $\zeta_\eps$ be the one parameter family of vector fields such
that $$ i_{\zeta_\eps} \om_\eps = \tau_\eps. $$ Now we integrate
the time dependent vector field $\zeta_\eps$ to a family
${\Psi_\eps}$ of diffeomorphisms (again we use that $M$ is
compact). We have

\begin{eqnarray*}
    \frac{d}{d\eps}\Psi^*_\eps \om_\eps & = & \Psi^*_\eps
\mathcal{L}_{\zeta_\eps} \om_\eps + \Psi_\eps^*\frac{d\om_\eps}{d\eps} \\
    & = & \Psi^*_\eps [d\iota_{\zeta_\eps} \om_\eps - d\tau_\eps] \\
    & = & \Psi^*_\eps [d\tau_\eps - d\tau_\eps] \\
    & = & 0.
\end{eqnarray*}

Hence $\Psi^*_\eps \om_\eps$ is constant so $\Psi^*_\eps \om_\eps = \om_0$ and
$\om_\eps = (\Psi_\eps^{-1})^* \om_0$.
\qed

\subsection{The geometric meaning of the family $\M$ of symplectic quotients -
shrinking the curvature}

We recall that $X_\ka$ denotes the complete simply-connected Riemannian manifold
of constant curvature $\ka$.  Let $r = (r_1, r_2, ... , r_n) \in (\R_+)^n$ with
$r$ not on a wall of $D_n$.   Let  $M_r(X_\ka)$ be the moduli space of n-gon
linkages with side-lengths $r$ in the space $X_\ka$.  The following theorem is
the main result of \cite{Sa}.

\begin{thm}
There exists $\al > 0$ and an analytically trivial fiber bundle $\pi:\e \to
(-\8,\al)$ such that $\pi^{-1}(\ka) = M_r(X_\ka)$.
\end{thm}

Let $\M$ be the family of symplectic quotients just constructed (except we will
take $(-\8,0]$ as base instead of $[1,0]$).  We then have

\begin{thm}
We have an isomorphism of fiber bundles
$$
\e|_{(-\8,0]} \simeq \M.
$$
\end{thm}
We will need
\begin{lem}
Let $\la>0$.  Then we have a canonical isomorphism
$$
M_r(X_\ka) \simeq M_{\la r}(X_{\ka/\la}).
$$
\end{lem}
\proof  Multiply the Riemannian metric on $X_\ka$ by $\la$.  Then the Riemannian
distance function is multiplied by $\la$ and the sectional curvature is
multiplied by $\frac{1}{\la}$.
\qed

\begin{rem}
There is a good way to visualize the above isomorphism by using the embedding of
$X_\ka, \;\ka<0$, in Minkowski space (as the upper sheet of the hyperboloid
$x^2+y^2+z^2-t^2 = \frac{-1}{\ka^2}$) or $X_\ka, \;\ka>0$ in $\R^4$ ( as the
sphere
$x^2+y^2+z^2+t^2 = \frac{1}{\ka^2}$).  The dilation map $v \mapsto \la v$ of the
ambient vector space maps $X_\ka$ to $X_{\ka/\la}$ and multiplies the
side-lengths
by $\la$.
\end{rem}

Now we can prove the theorem.  Let $\bar{p}: \M \to (-\8,0]$ be the family
constructed is \S\ref{17}.  By Lemma \ref{orbit} we see that $\bar{p}^{-1}(\eps)
\simeq
M_{\eps r}(X_{-1})$.  Thus we are shrinking the side-lengths of the n-gons as
$\eps \to 0$.  But we have just constructed a canonical isomorphism
$$
M_{\eps r}(X_{-1}) \simeq M_r(X_{-\eps}).
$$
So we may regard the deformation of \S\ref{17} as keeping the side-lengths fixed
and
shrinking the curvature to zero.

To give a formal proof we will construct an explicit diffeomorphism
$$
\xymatrix{
\M \ar[r] \ar[dr] & \e|_{(-\8,0]} \ar[d]\\
& (-\8,0]   }
$$

To this end, observe that the map $B\times (-\8,0] \to G/K \times (-\8,0]$ given
by $(b,\ka) \mapsto (b*K,\ka)$ induces a $K$-equivariant diffeomorphism

$$
\xymatrix{
B^n \times (-\8,0] \ar[r]^{F} \ar[dr]_{p} & (G/K)^n\times (-\8,0] \ar[d]^{\pi}\\
& (-\8,0]   }
$$
given by $F((b_1,...,b_n),\ka) = ((K, b_1*K, ..., b_1 * \cdots * b_{n-1}*K),
\ka)$.

We give $\pi^{-1}(\ka)$ the Riemannian metric $|\ka|(,)_\ka$.  Let
$\t{\e}' \subset (G/K)^n \times (-\8,0]$ be defined by $$ \t{\e}'
= \{(y_1, ... , y_n, \ka): \; y_1 = x_0, \; d_\ka(y_i,y_{i+1}) =
r_i, \; 1 \leq i \leq n \} $$ Here $d_\ka$ is the distance
function on $\pi^{-1}(\ka)$ associated to the Riemannian metric
$|\ka|(,)_\ka$.  Let $\Sigma_r \times (-\8,0]$ be the dressing
orbit through $(e^{r_1 H}, ... , e^{r_n H})$ for the $K$-dressing
action of $K$ on $p^{-1}(\ka)$.  We let $\t{\M} \subset \Sigma_r
\times (-\8,0]$ be the subset $\t{M} = \{(b_1,...,b_n, \ka): b_1*
b_2 * \cdots * b_n = 1\}$.  Then $F$ carries $\t{\M}$
diffeomorphically onto $\t{\e}'$ and induces the required
diffeomorphism $\M \to \e|_{(-\8,0]}$. \qed

\begin{rem}
The relative 2-form $\bar{\om}_\eps$ is a symplectic form along the fibers of
$\bar{p}$.
Thus
we have made the restriction of the family of \cite{Sa} to $(-\8,0]$ into a
family
of symplectic manifolds.  Can $\om_\eps$ be extended to $(-\8, \al)$for some
$\al>0$?
\end{rem}

\end{document}